%% file: arxiv_main.tex
\def\arxiv{}

\documentclass[mnsc,nonblindrev]{informs3} 

\OneAndAHalfSpacedXI 

\usepackage{endnotes}

 \usepackage{algpseudocode}
\usepackage{tikz}

\usepackage{natbib}
 \bibpunct[, ]{(}{)}{,}{a}{}{,}%
 \def\bibfont{\small}%

\usepackage[colorlinks=true,breaklinks=true,bookmarks=true,urlcolor=magenta,citecolor=magenta,linkcolor=magenta,bookmarksopen=false,draft=false]{hyperref}

\input{preamble}

\renewenvironment{proof}[1][Proof]{%
  \par\noindent\hspace{2em}{\itshape #1: }%
}{%
  \hspace*{\fill}~$\square$\par
}
\makeatother

\TheoremsNumberedThrough     
\ECRepeatTheorems

\EquationsNumberedThrough    

\begin{document}


\RUNAUTHOR{Jadbabaie, Shah, and Sinclair}

\RUNTITLE{Multi-Objective LQR with Linear Scalarization}

\TITLE{Multi-Objective LQR with Linear Scalarization}

\ARTICLEAUTHORS{%
\AUTHOR{Ali Jadbabaie, Devavrat Shah}
\AFF{
Institute for Data, Systems, and Society,
Massachusetts Institute of Technology, \\
Cambridge, MA 02139, \EMAIL{\{jadbabaie, devavrat\}@mit.edu}}
\AUTHOR{Sean R. Sinclair}
\AFF{Department of Industrial Engineering and Management Sciences, 
Northwestern University, \\
Evanston, IL 60208, \EMAIL{sean.sinclair@northwestern.edu}}
} 

\ABSTRACT{%
\input{r1_parts/abstract}
}%





\maketitle


\input{r1_parts/introduction}

\input{r1_parts/literature_review}
\input{r1_parts/preliminary}
\input{r1_parts/main_results}
\input{r1_parts/dare_sensitivity}
\input{r1_parts/main_proofs}
\input{r1_parts/conclusion}

\medskip

\noindent\textbf{Acknowledgments.} {Part of this work was done while Sean Sinclair was a Postdoctoral Associate at MIT under Devavrat Shah
and Ali Jadbabaie.  The authors would like to thank Swati Padmanabhan, Haoyuan Sun, and Daniel Pfrommer for insightful conversations about this work.}

\bibliographystyle{informs2014} 
\bibliography{references} 


%
%
%

\newpage

\crefalias{section}{appendix}
\begin{APPENDICES}
\OneAndAHalfSpacedXI 
    \input{r1_parts/appendix_lemmas}

\end{APPENDICES}








\end{document}

%% file: preamble.tex
\usepackage{amsmath,amsfonts,cancel}
\usepackage{thmtools}

\IfPackageLoadedTF{enumitem}{}{
\usepackage{enumitem}}
\usepackage{mathtools}
\usepackage{dsfont}
\usepackage{tcolorbox}
\usepackage{bbm}
\usepackage{xspace}
\usepackage{tensor}
\usepackage{bm}
\usepackage{tikz}
\usepackage{placeins}

\usepackage[capitalize]{cleveref}

\usepackage{color-edits}
\addauthor{srs}{orange}
\addauthor{TODO}{red}

\newcommand{\RR}{\mathbb{R}}
\renewcommand{\S}{\mathcal{S}}

\newcommand{\Exp}[1]{\mathbb E \left[ #1 \right]} 

\newcommand{\inv}{{-1}}
\newcommand{\T}{\mathcal{T}}

\renewcommand{\H}{\mathcal{H}}
\newcommand{\dare}{\textsf{dare}}
\newcommand{\dlyape}{\textsf{dlyape}}
\newcommand{\loss}{\mathcal{L}}
\newcommand{\hatloss}{\widehat{\mathcal{L}}}
\newcommand{\PF}{\textsf{PF}}
\newcommand{\hatPF}{\widehat{\textsf{PF}}}

\newcommand{\G}{\mathcal{G}}
\DeclareMathOperator{\Tr}{Tr}
\newcommand{\smin}{\sigma_{\min}}
\newcommand{\MObjLQR}{\textsf{MObjLQR}\xspace}
\newcommand{\smax}{\sigma_{\max}}

\renewcommand{\r}{R_{\max}}
\newcommand{\Pmax}{P_{\max}}
\newcommand{\Kmax}{K_{\max}}
\newcommand{\Khat}{\widehat{K}}
\newcommand{\Ahat}{\widehat{A}\xspace}
\newcommand{\Bhat}{\widehat{B}\xspace}

\newcommand{\gbar}{\overline{\gamma}}
\newcommand{\tbar}{\overline{\tau}}

\usepackage{multirow}

\crefname{assumption}{assumption}{assumptions}
\crefname{algocf}{algorithm}{algorithm}

\mathchardef\mhyphen="2D 
\ifdefined \argmin
\else \DeclareMathOperator*{\argmin}{arg\,min}
\fi
\ifdefined \argmax
\else \DeclareMathOperator*{\argmax}{arg\,max}
\fi
\ifdefined \norm \else \DeclarePairedDelimiter{\norm}{\lVert}{\rVert} \fi
\ifdefined \ceil \else \DeclarePairedDelimiter{\ceil}{\lceil}{\rceil} \fi

\let\originalleft\left
\let\originalright\right
\renewcommand{\left}{\mathopen{}\mathclose\bgroup\originalleft}
\renewcommand{\right}{\aftergroup\egroup\originalright}



%% file: r1_parts/abstract.tex
We study finite approximation of the Pareto front for the multi-objective linear quadratic regulator (LQR). Classical results characterize Pareto-optimal controllers through weighted sums of the objectives, but do not quantify the effect of discretizing the scalarization weights. We prove the needed sensitivity property, showing that an $\epsilon$ perturbation of the weights changes the full objective vector by $O(\epsilon)$. Consequently, computing optimal controls over a finite grid of the weight simplex uniformly approximates the Pareto front. We also extend the guarantee to unknown dynamics via certainty equivalence, showing that if the model estimates are accurate to the target accuracy, the same procedure stabilizes the true system and attains the same Pareto-front guarantee.

%% file: r1_parts/introduction.tex
\section{Introduction}
\label{sec:intro}

Modern control and optimization techniques shape much of our physical and digital infrastructure, including data centers, power grids, and supply chains. Many of these systems involve multiple objectives chosen by stakeholders. For example, control methods have been used to manage power distribution in plug-in electric vehicles by optimizing a handcrafted objective to improve system efficiency and vehicle response~\citep{lu2019multi}. Similar tradeoffs arise in dynamic asset allocation, where one balances risk and return~\citep{sridharan2011receding}, and in smart-grid energy management, where algorithms trade off energy use, cost, reliability, and environmental impact~\citep{shaikh2016intelligent,yang2012multi}.

Traditional approaches for multiple objectives rely on heuristics or stakeholder input to construct a single objective. This can obscure tradeoffs between metrics and does not reflect that there is typically no single best ranking of solutions. Instead, a preference-agnostic approach returns a {\em set} of solutions. We focus on Pareto optimality (see \cref{def:pareto_front}), where a solution is Pareto optimal if no other solution improves one objective without worsening another. This creates computational challenges because the Pareto front is continuous, so a discrete approximation requires smoothness~\citep{das1997closer}. 
%

\subsection{Our Contributions}

Recent work has studied multi-objective reinforcement learning in finite-horizon or discounted infinite-horizon settings~\citep{liu2014multiobjective}. Less is known in undiscounted infinite-horizon settings, where stability is central. We study the Linear Quadratic Regulator (LQR), a canonical model for this regime~\citep{anderson2007optimal}. In multi-objective LQR, prior work shows that linear scalarization characterizes the Pareto front~\citep{khargonekar2002multiple,rotea1991h2}. This suggests a simple algorithm: discretize the scalarization weights and solve the single-objective LQR problem per grid point. This approach has been used in practice without theoretical guarantees~\citep{logist2009efficient}. Our main contributions are as follows:

\smallskip \noindent {\bf Approximation guarantee} (\cref{thm:pf_approx}). We prove that the discretization algorithm uniformly approximates the Pareto front by establishing smoothness with respect to the scalarization weights. Our analysis builds on Riccati sensitivity bounds from \citet{mania2019certainty}, which we extend to perturbations in the cost matrices. See \cref{fig:sample} for an illustration on a multi-objective pendulum problem.

\smallskip \noindent {\bf Certainty equivalence under unknown dynamics} (\cref{thm:ce_pf_approximation}). We extend the approximation guarantee to settings where the system dynamics are unknown and replaced by estimates. We show that the same algorithm, applied with certainty-equivalent dynamics, approximates the Pareto front when the estimation error is of the same order as the desired approximation accuracy.

\ifdefined\arxiv
The techniques combine ideas from semidefinite programming, optimal control, and mathematical optimization. They provide a theoretical justification for a simple Pareto-front approximation method already used in practice~\citep{logist2009efficient}, while relying only on standard single-objective LQR solvers~\citep{alessio2009survey}. We expect these methods to be useful in settings with linear dynamics, including congestion control~\citep{bartoszewicz2013linear}, queueing systems~\citep{bonald2015multi}, and inventory control~\citep{parlar2003balancing}.
\else\fi

%% file: r1_parts/literature_review.tex
\subsection{Related Work}
\label{sec:related_work}

Our paper builds on work in multi-objective control, scalarization, and data-driven LQR. For background on LQR see \citet{anderson2007optimal}, and for Pareto optimality and scalarization see \citet{miettinen1999nonlinear}.

\smallskip
\noindent \textbf{Characterizing the Pareto frontier.} Pareto-front computation for control systems often uses scalarization and related methods~\citep{logist2009efficient,li1994multilevel,van2008efficiently,li1993general,koussoulas1986multiple,nair2024multi}. Chebyshev scalarization and the $\epsilon$-constraint method are common because they can recover non-convex Pareto fronts~\citep{miettinen1999nonlinear,das1997closer}. Linear scalarization and separating-hyperplane arguments date back at least to \citet{da1967constrained}, and scalarization characterizations for linear quadratic, linear-quadratic-Gaussian, and $H_2$-type control appear in classical work including \citet{koussoulas1986multiple,makila1989multiple,rotea1991h2,li1993general,khargonekar2002multiple}. We build on this characterization by showing that discretizing the scalarization simplex yields a uniform approximation of the Pareto front.

\smallskip
\noindent \textbf{Empirical applications of multi-objective control.} Linear scalarization is widely used in multi-objective control applications because it is simple and preserves tractable single-objective subproblems. Examples include vehicle suspension control~\citep{genov2019linear,ye2019comparative}, aircraft control~\citep{giacoman2008robust}, and switched systems~\citep{wu2018optimal}. Our results provide a theoretical guarantee for this practical procedure in multi-objective LQR.

\smallskip
\noindent \textbf{Data-driven control.} A separate line of work studies data-driven control with a single objective. Policy-gradient methods for LQR have been analyzed using global landscape properties of the LQR objective~\citep{fazel2018global}. Certainty-equivalent and adaptive LQR methods synthesize controllers from estimated dynamics using samples~\citep{arora2018towards,cohen2018online}, advice or predictions~\citep{yu2020power,zhang2021regret}, or adversarial models~\citep{agarwal2019online,li2021online}. Our results use sensitivity bounds for the discrete algebraic Riccati equation, whose perturbation theory has a long history~\citep{konstantinov1993perturbation,sun1998perturbation,sun1998sensitivity}. We extend the results in \citet{mania2019certainty} to account for perturbations in the cost matrices.

%% file: r1_parts/preliminary.tex
\begin{figure}[!t]
    \centering
    \ifdefined\arxiv
    \includegraphics[width=0.7\linewidth]{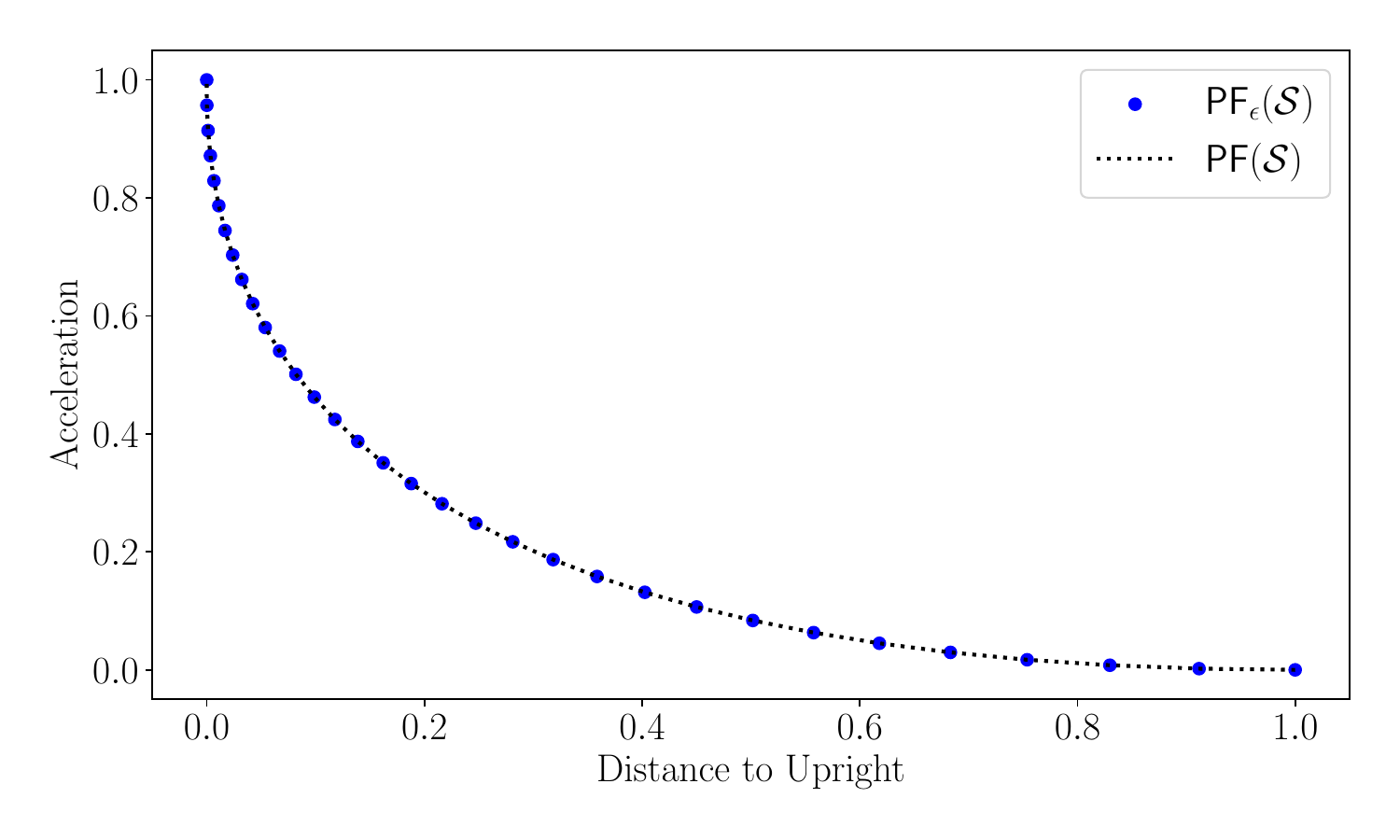}
    \else
    \includegraphics[width=\linewidth]{figures/pendulum_tradeoff.pdf}
    \fi
    \caption{
    Approximation of the Pareto front for an inverted pendulum. The axes show normalized performance on distance-to-upright and cumulative-acceleration objectives. Each point is a controller computed by the algorithm in \cref{sec:main_results} with known dynamics and $\epsilon=10^{-1.5}$.
    }
    \label{fig:sample}
\end{figure}

\section{Preliminaries and Notation}
\label{sec:preliminary}

In the following we write $[m]=\{1,\ldots,m\}$ and $\Delta([m])=\{w\in\mathbb{R}^m_{\ge0}:\norm{w}_1=1\}$. Let $\lambda_i(A)$ denote the eigenvalues of $A$, ordered by non-increasing absolute value, and $\sigma_i(A)$ its singular values, ordered non-increasingly. We write $\smin(A)$ and $\smax(A)$ for the minimum and maximum singular values, and $\rho(A)$ for the spectral radius. Unless stated otherwise, matrix norms are spectral operator norms, and $O(1)$ denotes absolute constants.  Lastly, we use $\norm{\cdot}_+ = 1 + \norm{\cdot}$.

\subsection{Multi-Objective LQR}
The standard linear quadratic regulator (LQR) is parameterized by dynamics matrices $A \in \RR^{n\times n}$ and $B \in \RR^{n\times d}$. We assume $n \ge d$ for simplicity, though our guarantees scale with $\min\{n,d\}$. Given state $x_t$ and control $u_t$, the state evolves as $x_{t+1} = A x_t + B u_t + \xi_t$ where $\xi_t \sim N(0, \sigma^2 I_n)$.
Classical LQR uses a single quadratic cost $(Q,R)$, whereas we consider multiple cost matrices $(Q_i,R_i)_{i\in[m]}$, with $Q_i\in\RR^{n\times n}$ and $R_i\in\RR^{d\times d}$. Since the optimal controller for any single-objective LQR problem is linear~\citep{anderson2007optimal}, we restrict attention to controls $u_t=Kx_t$, where $K\in\RR^{d\times n}$. For fixed cost matrices $Q$ and $R$, define
\begin{align}
\label{eq:lqr_loss}
\loss(K,Q,R) & \coloneq \lim_{T\to\infty} \mathbb{E} \left[ \frac{1}{T}\sum_{t=0}^T \left(x_t^\top Qx_t+u_t^\top Ru_t\right) \right] \\
\text{s.t.}\quad x_{t+1} & = Ax_t + Bu_t + \xi_t, \qquad u_t = Kx_t, \notag
\end{align}
where the expectation is over $x_0 \sim N(0,I_n)$ and the noise terms $(\xi_t)_{t \ge 0}$. We make the following assumption:

\begin{assumption}
\label{ass:stabilizable_observable}
The pair $(A,B)$ is stabilizable, and each $Q_i$ and $R_i$ is positive definite. Since rescaling $(Q_i,R_i)$ does not change the optimal controller, we assume without loss of generality that $\smin(Q_i),\smin(R_i)\ge1$ for all $i\in[m]$.
\end{assumption}

For the remainder of the paper, we assume \Cref{ass:stabilizable_observable} holds. This assumption implies detectability of $(A,Q_i^{1/2})$ for each $i\in[m]$~\citep{anderson2007optimal}. We define the set of stabilizing controls as
\begin{equation}
\label{eq:stable_set}
\S = \{K \mid \rho(A+BK)<1\}.
\end{equation}
The goal of {\bf multi-objective LQR} (\MObjLQR) is to solve
\begin{align}
\label{eq:mobj_lqr}
    \min_{K\in\S} \, \{ \loss(K,Q_i,R_i) \}_{i\in[m]} .
\end{align}
The ``optimal'' control, in this case, is unclear since the control which optimizes $(Q_j,R_j)$ need not perform well for $(Q_i,R_i)$ with $i\neq j$.  When $m=1$, or equivalently when the cost matrices $(Q,R)$ are fixed, the optimal policy is the linear feedback control $u_t=Kx_t$ with
\begin{equation}
\label{eq:optimal_solution}
K = - (R + B^\top P B)^{-1} B^\top P A,
\end{equation}
where $P$ is the positive definite solution of the discrete algebraic Riccati equation
\begin{equation}
\label{eq:dare}
P = A^\top P A - A^\top P B(R + B^\top P B)^{-1} B^\top P A + Q.
\end{equation}
We write $\dare(A,B,Q,R)$ for the unique positive semidefinite solution to \Cref{eq:dare}~\citep{anderson2007optimal}. We also assume $\smin(P)\ge1$ without loss of generality; this follows from \Cref{ass:stabilizable_observable}, \cref{lem:eigenvalue_dare_bound} (appendix), and \citet{garloff1986bounds}.

\subsection{Pareto Optimality and Linear Scalarization}
Our goal is to find controls $K\in\S$ with strong performance {\em uniformly} across the objectives $\loss_i(K)\coloneq \loss(K,Q_i,R_i)$ for $i\in[m]$. We consider {\em Pareto optimality} as our solution concept:

\begin{definition}[Pareto optimal]
\label{def:pareto_front}
A control $K\in\S$ is {\bf Pareto optimal} if there is no $K'\in\S$ such that $\loss_i(K')\le \loss_i(K)$ for all $i\in[m]$ and $\loss_j(K')<\loss_j(K)$ for some $j\in[m]$. We denote the set of Pareto optimal controls by
\[
    \PF(\S)=\{K\in\S: K \text{ is Pareto optimal}\}.
\]
\end{definition}

A standard way to generate Pareto optimal solutions is linear scalarization. For a weight $w\in\Delta([m])$, define
\begin{align}
\label{eq:linear_scalarize}
    \loss_w(K)\coloneq \sum_i w_i\loss_i(K).
\end{align}

The key feature of \MObjLQR is that this scalarized objective is itself an LQR objective with weighted cost matrices.

\begin{lemma}
\label{lem:linearity_lqr_objective}
For any stabilizing control $K\in\S$,
\[
    \loss_w(K)
    =
    \sum_i w_i\loss_i(K)
    =
    \loss(K,Q_w,R_w),
\]
where $Q_w = \sum_i w_i Q_i$ and $R_w = \sum_i w_i R_i$.
\end{lemma}
\begin{proof}
The result follows by linearity of the stage cost:
\[
\sum_i w_i(x_t^\top Q_i x_t+u_t^\top R_i u_t)
=
x_t^\top Q_wx_t+u_t^\top R_wu_t.
\]
Since $K\in\S$, the average costs are finite, so the finite sum over objectives can be interchanged with the limit.
\end{proof}

Thus, for each fixed $w$, optimizing $\loss_w$ reduces to a standard single-objective LQR problem with cost matrices $(Q_w,R_w)$. If $K_w$ denotes this optimizer, then
\[
    K_w=-(R_w + B^\top P_w B)^{-1} B^\top P_w A,
\]
where $P_w=\dare(A,B,Q_w,R_w)$. Hence linear scalarization preserves the Riccati structure of LQR and allows one to use existing numerical solvers~\citep{alessio2009survey}.
We lastly recall a classical result from multi-objective control that the covariance/semidefinite formulation convexifies the attainable objective set, and linear scalarization therefore recovers the Pareto-optimal stabilizing controllers; see
\citet{khargonekar2002multiple,rotea1991h2,makila1989multiple}.

\begin{theorem}
\label{thm:sufficiency_linear}
A stabilizing controller $K\in\S$ is Pareto optimal if and only if there exists
$w\in\Delta([m])$ such that $K$ minimizes $\loss(K,Q_w,R_w)$ over $\S$. Equivalently,
\[
    \PF(\S)=\{K_w : w\in\Delta([m])\}.
\]
\end{theorem}

Thus, approximating the Pareto front reduces to understanding how the optimal controller $K_w$ varies with the scalarization weight $w$. We address this question in the next section.

%% file: r1_parts/main_results.tex
\section{Main Results}
\label{sec:main_results}

The scalarization result in \cref{thm:sufficiency_linear} reduces Pareto-front approximation to a one-parameter family of standard LQR problems.  For each weight $w\in\Delta([m])$, let $P_w = \dare(A,B,Q_w, R_w)$ and
\begin{align*}
K_w & =-(R_w+B^\top P_wB)^{-1}B^\top P_wA .
\end{align*}
Then every Pareto-optimal controller can be represented as $K_w$ for some $w$.  This suggests the following approach.

\smallskip
\noindent\textbf{Algorithm.}
Let $N_\epsilon$ be an $\epsilon$-net of $\Delta([m])$ in $\ell_1$, so that for every $w\in\Delta([m])$ there exists $w_\epsilon\in N_\epsilon$ with $\norm{w-w_\epsilon}_1\leq\epsilon$.  Such a net can be chosen with $|N_\epsilon|=O(\epsilon^{-m})$.  For each $w_\epsilon\in N_\epsilon$, solve the single-objective LQR problem with cost matrices $Q_{w_\epsilon}$ and $R_{w_\epsilon}$, and
return
\begin{equation}
\label{eq:approx_ccs}
\PF_\epsilon(\S)
=
\{K_{w_\epsilon}:w_\epsilon\in N_\epsilon\}.
\end{equation}
The algorithm requires $O(\epsilon^{-m})$ calls to a standard LQR solver.
Discretizing weights alone, however, does not guarantee objective-space approximation. In general multi-objective problems, nearby scalarization weights can produce distant Pareto points~\citep{das1997closer}. 

\smallskip
\noindent\textbf{Uniform constants.}
Before presenting our main result, we need constants that control all scalarized LQR instances simultaneously, rather than constants that depend on a particular weight $w$. We use the following quantities:
\[
\Pmax=\sup_{w\in\Delta([m])}\norm{P_w},
\qquad
\Kmax=\sup_{w\in\Delta([m])}\norm{K_w},
\]
and collect the primitive bounds into
\ifdefined\arxiv
\begin{align*}
\Gamma
=
\max\{
1+\Pmax,\,
1+\Kmax,\,
\norm{A}_+,\,
\norm{B}_+, 
\max_{i\in[m]}\norm{Q_i},\,
\max_{i\in[m]}\norm{R_i},\,
1+m\max_{i\in[m]}\norm{R_i^{-1}}\}.
\end{align*}
\else
\begin{align*}
\Gamma
=
\max\{&
1+\Pmax,\,
1+\Kmax,\,
\norm{A}_+,\,
\norm{B}_+, \\
&
\max_{i\in[m]}\norm{Q_i},\,
\max_{i\in[m]}\norm{R_i},\,
1+m\max_{i\in[m]}\norm{R_i^{-1}}\}.
\end{align*}
\fi
Note that $\Gamma$ has no dependence on a particular $w \in \Delta([m])$.  The quantities $\Pmax$ and $\Kmax$ are finite, and \cref{lem:upper_bounds} (\cref{sec:dare_constants}) gives explicit upper bounds.  We also need a uniform closed-loop stability margin over the scalarized optimizers.  In \cref{lem:uniform_control_pf} (\cref{sec:dare_constants}) , we show that there exist constants $\gbar<1$ and $\tbar<\infty$ such that
\[
\rho(A+BK_w)\leq \gbar,
\qquad
\tau(A+BK_w,\gbar)\leq \tbar
\]
for all $w\in\Delta([m])$.

Our main theorem shows that the finite set in \cref{eq:approx_ccs} uniformly approximates the Pareto front in the original objective coordinates.
\begin{theorem}
\label{thm:pf_approx}
Suppose that \Cref{ass:stabilizable_observable} holds and
\[
\epsilon
\leq
O(1)
\frac{(1-\gbar^2)^4}{\tbar^4}
\Gamma^{-17}.
\]
Then, for every $K\in\PF(\S)$, there exists
$K_\epsilon\in\PF_\epsilon(\S)$ with
\[
\norm{\vec{\loss}(K)-\vec{\loss}(K_\epsilon)}_\infty
\leq
O(1)
\sigma^2 n
\Gamma^{29}
\frac{\tbar^8}{(1-\gbar^2)^4}
\epsilon .
\]
\end{theorem}

This guarantee controls each original objective separately.  It shows that for any Pareto-optimal controller $K_w$, solving the LQR problem at a nearby grid point $w_\epsilon$ produces a controller $K_{w_\epsilon}$ whose full loss vector is close to that of $K_w$.  The polynomial dependence on $\Gamma$, $\tbar$, and $(1-\gbar^2)^{-1}$ is unlikely to be tight, but it makes explicit that the guarantee worsens as the closed-loop systems approach instability~\citep{mania2019certainty}. The proof is deferred to \cref{sec:main_proofs}.

\smallskip
\noindent\textbf{Certainty-equivalent approximation.}
We next consider the case where the dynamics are unknown.  Suppose the algorithm has access to estimates $\Ahat$ and $\Bhat$ of $A$ and $B$.  For fixed cost matrices $Q$ and $R$, define the estimated-dynamics loss
\begin{align}
\label{eq:approx_LQR_loss}
\hatloss(K,Q,R)
\coloneq
\lim_{T\rightarrow\infty}
\mathbb{E}
\left[
\frac{1}{T}
\sum_{t=0}^{T}
\left(
x_t^\top Qx_t+u_t^\top Ru_t
\right)
\right], \\
\text{s.t.}\quad x_{t+1}=\Ahat x_t+\Bhat u_t+\xi_t,
\qquad
u_t=Kx_t. \notag
\end{align}
\noindent Write
$\hatloss_i(K)=\hatloss(K,Q_i,R_i)$ and
$\hatloss_w(K)=\hatloss(K,Q_w,R_w)$.

The certainty-equivalent algorithm applies the same simplex discretization, but treats the estimated dynamics as exact.  For each $w_\epsilon\in N_\epsilon$, compute
\[
\Khat_{w_\epsilon}
\in
\argmin_K \hatloss_{w_\epsilon}(K)
\]
and return
\begin{equation}
\label{eq:approx_ccs_ce}
\hatPF_\epsilon(\S)
=
\{\Khat_{w_\epsilon}:w_\epsilon\in N_\epsilon\}.
\end{equation}
The controller $\Khat_{w_\epsilon}$ is constructed to stabilize the estimated dynamics.  The next result shows that, when the estimation error is of the same order as the target approximation accuracy, it also stabilizes the true dynamics and preserves the Pareto-front approximation guarantee.

\begin{theorem}
\label{thm:ce_pf_approximation}
Suppose that $\norm{\Ahat - A} \leq \epsilon$ and $\norm{\Bhat - B} \leq \epsilon$, $(\Ahat,\Bhat)$ is stabilizable, \Cref{ass:stabilizable_observable} holds, and
\[
\epsilon
\leq
O(1)
\frac{(1-\gbar^2)^4}{\tbar^4}
\Gamma^{-17}.
\]
Then, for every $K\in\PF(\S)$, there exists $\Khat_\epsilon\in\hatPF_\epsilon(\S)$ such that $\Khat_\epsilon$ stabilizes the true dynamics $(A,B)$ and
\[
\norm{
\vec{\loss}(K)-\vec{\loss}(\Khat_\epsilon)
}_\infty
\leq
O(1)
\sigma^2 n
\Gamma^{29}
\frac{\tbar^8}{(1-\gbar^2)^4}
\epsilon .
\]
\end{theorem}

\begin{proof}[Proof Sketch]%
The proofs of \cref{thm:pf_approx,thm:ce_pf_approximation} have three ingredients.  First, \cref{sec:dare_sensitivity} establishes perturbation bounds for the discrete algebraic Riccati equation when $(A,B,Q,R)$ are perturbed.  Second, \cref{sec:dare_constants} converts those local perturbation bounds into $w$-independent bounds over the scalarization simplex using $\Gamma$, $\gbar$, and $\tbar$.  Third, the proofs in \cref{sec:main_proofs} convert Riccati sensitivity into controller sensitivity and then into objective-wise loss bounds.  The certainty-equivalent theorem follows from the same perturbation argument, with the additional observation that sufficiently small model error preserves stability under the true dynamics.
\end{proof}

%% file: r1_parts/dare_sensitivity.tex
\section{Perturbation Theory for the Discrete Riccati Equation}
\label{sec:dare_sensitivity}

A necessary aspect of our analysis is perturbation theory for solutions to the discrete algebraic Riccati equation as the problem parameters vary. More specifically, let $P=\dare(A,B,Q,R)$ and $P_\epsilon=\dare(A_\epsilon,B_\epsilon,Q_\epsilon,R_\epsilon)$, where each pair of matrices differs by at most $\epsilon$ (i.e., $\norm{A-A_\epsilon}\le\epsilon$, etc.). Our goal is to show that $\norm{P_\epsilon-P}\lesssim\epsilon$ for sufficiently small $\epsilon$.  There is a large body of literature on perturbation theory for algebraic Riccati equations~\citep{sun1998perturbation,sun1998sensitivity}. We adapt the operator-theoretic approach of~\citet{mania2019certainty,konstantinov1993perturbation} to account for perturbations in the cost matrices $(Q,R)$. This argument yields explicit constants that can be tracked uniformly over the scalarization simplex and propagated to the final approximation bound.

\smallskip
\noindent \textbf{Notation.} We use $\Delta_M=M_\epsilon-M$ to denote matrix perturbations, where $M$ can vary. We also denote $\norm{\cdot}_+=\norm{\cdot}+1$, and let $\r=1+\max\{\norm{R^\inv},\norm{R_\epsilon^\inv}\}$. For a matrix $L$, define
\begin{equation}
\label{eq:tau_def}
\tau(L,\rho)=\sup_{k\ge0}\, \norm{L^k}\rho^{-k}
\end{equation}
as its rate of growth.  Equivalently, $\tau(L,\rho)$ is the smallest value such that $\norm{L^k}\le\tau(L,\rho)\rho^k$ for all $k\ge0$. This quantity may be infinite, but by Gelfand's formula it is finite whenever $\rho>\rho(L)$~\citep{anderson2007optimal}. Thus, if $L$ is stable, one can choose $\rho<1$ with $\tau(L,\rho)<\infty$. At a high level, $\tau(L,\rho)$ quantifies how quickly the closed-loop system contracts; less stable systems have larger values. See \citet{tu2017non,mania2019certainty} for more details. Our primary goal in this section is to establish the following:

\begin{theorem}
\label{thm:dare_sensitivity}
Let $P = \dare(A,B,Q,R)$, $P_\epsilon = \dare(A_\epsilon,B_\epsilon,Q_\epsilon, R_\epsilon)$, and let
\[
K=-(R+B^\top PB)^{-1}B^\top PA,
\qquad
L=A+BK .
\]
We assume that $(A,B)$ and $(A_\epsilon,B_\epsilon)$ are stabilizable, $Q$ and $Q_\epsilon$ are positive definite, $R$ and $R_\epsilon$ are positive definite, $\smin(P) \geq 1$, and $\max\{\norm{\Delta_A}, \norm{\Delta_B}, \norm{\Delta_Q}, \norm{\Delta_R}\} \leq \epsilon$.  Then:
\[
\norm{P_\epsilon - P} \leq O(1) \frac{\tau(L, \rho)^2}{1 - \rho^2} \norm{A}_+^2 \norm{P}_+^2 \norm{B}_+^2 \r^2 \epsilon,
\]
so long as
\[
\epsilon \leq O(1) \frac{(1 - \rho^2)^4}{\tau(L, \rho)^4} \norm{A}_+^{-2} \norm{P}_+^{-3} \norm{B}_+^{-4} \r^{-3} \norm{L}_+^{-2}.
\]
\end{theorem}

Although \cref{thm:dare_sensitivity} gives the desired $O(\epsilon)$ perturbation bound, its constants still depend on the particular instance $(A,B,Q,R)$. In our approximation algorithm, this instance varies with the scalarization weight.  For $w\in\Delta([m])$ and its nearest grid point $w_\epsilon\in N_\epsilon$, we apply the theorem with $P=\dare(A,B,Q_w,R_w)$ and $P_\epsilon=\dare(A,B,Q_{w_\epsilon},R_{w_\epsilon})$.
Thus the resulting constants may depend on $w$, whereas our final approximation guarantee should depend only on the primitives $(A,B,(Q_i,R_i)_{i\in[m]})$. We address this uniformity issue in \cref{sec:dare_constants}. Before proving \cref{thm:dare_sensitivity} in \cref{sec:dare_proof}, we give a brief proof sketch.

\begin{proof}[Proof Sketch]
Denote by $F(X,A,B,Q,R)$ the matrix expression
\ifdefined\arxiv
\begin{align}
\label{eq:f_definition}
F(X,A,B,Q,R)
&= X - A^\top X A 
+ A^\top X B(R+B^\top X B)^{-1}B^\top X A
- Q \notag \\
&= X - A^\top X(I+BR^{-1}B^\top X)^{-1}A - Q . 
\end{align}
\else
\begin{align}
\label{eq:f_definition}
F(X,A,B,Q,R)
&= X - A^\top X A \\
&\quad
+ A^\top X B(R+B^\top X B)^{-1}B^\top X A
- Q \notag \\
&= X - A^\top X(I+BR^{-1}B^\top X)^{-1}A - Q . \notag
\end{align}
\fi
Solving the Riccati equation for $(A,B,Q,R)$ is equivalent to finding the unique positive definite $X$ with $F(X,A,B,Q,R)=0$. Let $\Delta_P=P_\epsilon-P$. Since $P$ and $P_\epsilon$ solve their respective DAREs, $F(P,A,B,Q,R)=0$ and $F(P_\epsilon,A_\epsilon,B_\epsilon,Q_\epsilon,R_\epsilon)=0$.

We will start by constructing an operator $\Phi(X)$ such that any fixed point $X$ satisfying $X = \Phi(X)$ must be equal to $\Delta_P$.  However, to show that $\Phi$ has a fixed point, we consider a set $S_\nu$ of matrices satisfying $\norm{X} \leq \nu$.  We will show that $\Phi$ maps $S_\nu$ into itself, and is a contraction.  Hence, $\Phi$ has a fixed point equal to $\Delta_P$, and since $\Delta_P \in S_\nu$ must satisfy $\norm{\Delta_P} \leq \nu$.  The proof finishes by selecting $\nu = O(\epsilon)$.
\end{proof}

\subsection{Proof of Riccati Sensitivity}
\label{sec:dare_proof}

\begin{proof}[Proof of \cref{thm:dare_sensitivity}]
For convenience we use $S = B R^\inv B^\top$ and $S_\epsilon = B_\epsilon R_\epsilon^\inv B_\epsilon^\top$.  For any matrix $X$ such that $I+S(P+X)$ is invertible, note that
\ifdefined\arxiv
\begin{align}
\label{eq:f_inv_x}
F(P+X,A,B,Q,R) & = X - L^\top X L + L^\top X(I+S(P+X))^\inv SXL,
\end{align}
\else
\begin{align}
\label{eq:f_inv_x}
F(P+X,A,B,Q,R) & = X - L^\top X L \\
&\quad + L^\top X(I+S(P+X))^\inv SXL, \notag
\end{align}
\fi
where $L=A+BK$.  This follows from adding $F(P,A,B,Q,R)$, which is equal to zero, to the right-hand side of \cref{eq:f_inv_x} and using that $(I+BR^\inv B^\top P)^\inv A = A + BK$.
Denote
\begin{align*}
\T(X)=X-L^\top XL, \,\,\, \H(X)=L^\top X(I+S(P+X))^\inv SXL .
\end{align*}
Then \cref{eq:f_inv_x} says that
\[
F(P+X,A,B,Q,R)=\T(X)+\H(X).
\]
Since \cref{eq:f_inv_x} is satisfied for any matrix $X$ such that $I+S(P+X)$ is invertible, the matrix equation
\ifdefined\arxiv
\begin{align}
\label{eq:f_difference}
F(P+X,A,B,Q,R)
-F(P+X,A_\epsilon,B_\epsilon,Q_\epsilon,R_\epsilon)
=
\T(X)+\H(X)
\end{align}
\else
\begin{align}
\label{eq:f_difference}
&F(P+X,A,B,Q,R) \notag \\
&\quad
-F(P+X,A_\epsilon,B_\epsilon,Q_\epsilon,R_\epsilon)
=
\T(X)+\H(X)
\end{align}
\fi
has a unique symmetric solution $X$ such that $P+X\succeq0$.  This solution is $X=\Delta_P$ because any solution must satisfy that $F(P+X,A_\epsilon,B_\epsilon,Q_\epsilon,R_\epsilon)=0$.

Note that the linear map $\T : X \rightarrow X - L^\top X L$ has eigenvalues equal to $1-\lambda_i\lambda_j$, where $\lambda_i$ and $\lambda_j$ are eigenvalues of the matrix $L$.  Since $L=A+BK$ is stable, the linear map $\T$ is invertible.  We define
\ifdefined\arxiv
\begin{align*}
\Phi(X)
=
\T^\inv\Big(
F(P+X,A,B,Q,R) 
-
F(P+X,A_\epsilon,B_\epsilon,Q_\epsilon,R_\epsilon)
-\H(X)
\Big).
\end{align*}

\else
\begin{align*}
\Phi(X)
&=
\T^\inv\Big(
F(P+X,A,B,Q,R) \\
&\qquad
-
F(P+X,A_\epsilon,B_\epsilon,Q_\epsilon,R_\epsilon)
-\H(X)
\Big).
\end{align*}
\fi
Then solving for $X$ in \cref{eq:f_difference} is equivalent to finding an $X$ satisfying $P+X\succeq0$ with $X=\Phi(X)$.  Thus we have that $\Phi$ has a unique symmetric fixed point $X$ such that $P+X\succeq0$ and that is $X=\Delta_P$.

The remainder of the proof is focused on establishing that $\Phi(X)$ has a fixed point with bounded norm. Define the set
\[
\S_\nu = \{X\mid \norm{X}\leq\nu,\ X=X^\top,\ P+X\succeq0\}.
\]
By assumption we have that $\norm{\Delta_A},\norm{\Delta_B},\norm{\Delta_Q},\norm{\Delta_R}\leq\epsilon$. Hence, we also have that
\ifdefined\arxiv
\begin{align*}
\norm{\Delta_S} & \leq \norm{B_\epsilon R^\inv B_\epsilon^\top-BR^\inv B^\top} 
+\norm{B_\epsilon R_\epsilon^\inv B_\epsilon^\top
-B_\epsilon R^\inv B_\epsilon^\top} \\
&\leq
\norm{\Delta_B R^\inv\Delta_B^\top}
+\norm{B R^\inv\Delta_B^\top} 
+\norm{\Delta_B R^\inv B^\top}
+\norm{B_\epsilon}^2\norm{R_\epsilon^\inv-R^\inv} \\
&\leq
\epsilon^2\norm{R^\inv}
+2\norm{B}\norm{R^\inv}\epsilon 
+(1+\norm{B})^2\norm{R_\epsilon^\inv-R^\inv} \\
&\leq
3(1+\norm{B})\norm{R^\inv}\epsilon 
+(1+\norm{B})^2\norm{R_\epsilon^\inv-R^\inv},
\end{align*}
\else
\begin{align*}
\norm{\Delta_S} & \leq \norm{B_\epsilon R^\inv B_\epsilon^\top-BR^\inv B^\top} \\
&\quad
+\norm{B_\epsilon R_\epsilon^\inv B_\epsilon^\top
-B_\epsilon R^\inv B_\epsilon^\top} \\
&\leq
\norm{\Delta_B R^\inv\Delta_B^\top}
+\norm{B R^\inv\Delta_B^\top} \\
&\quad
+\norm{\Delta_B R^\inv B^\top}
+\norm{B_\epsilon}^2\norm{R_\epsilon^\inv-R^\inv} \\
&\leq
\epsilon^2\norm{R^\inv}
+2\norm{B}\norm{R^\inv}\epsilon \\
&\quad
+(1+\norm{B})^2\norm{R_\epsilon^\inv-R^\inv} \\
&\leq
3(1+\norm{B})\norm{R^\inv}\epsilon \\
&\quad
+(1+\norm{B})^2\norm{R_\epsilon^\inv-R^\inv},
\end{align*}
\fi
where in the last line we used $\epsilon \leq 1$.  However,
\begin{align*}
\norm{R_\epsilon^\inv-R^\inv}
&=
\norm{R_\epsilon^\inv(R-R_\epsilon)R^\inv} \\
&\leq
\norm{R^\inv}\norm{R_\epsilon^\inv}\norm{R-R_\epsilon}
\leq
\r^2\epsilon .
\end{align*}
Combining this with before, we have
\(
\norm{\Delta_S}\leq4\norm{B}_+^2\r^2\epsilon.
\)

With this in hand, we next show that $\Phi$ maps $\S_\nu$ to itself, and is a contraction.  The Banach fixed point theorem then implies that $\Phi$ has a fixed point, and that fixed point must belong to $\S_\nu$.

\begin{lemma}
\label{lem:contraction}
Suppose that $X,X_1,$ and $X_2$ are in $\S_\nu$ for $\nu \leq \min\{1,\norm{S}^\inv\}$.  Moreover, suppose $\epsilon \leq 1$ and $\sigma_{min}(P) \geq 1$.  Then:
\begin{align*}
\norm{\Phi(X)}
&\leq C_1\nu^2+C_2\epsilon, \\
\norm{\Phi(X_1)-\Phi(X_2)}
&\leq
\norm{X_1-X_2}(C_3\nu+C_4\epsilon),
\end{align*}
where
$C_1=\frac{\tau(L,\rho)^2}{1-\rho^2}\norm{L}^2\norm{S}$,
$C_2=8\frac{\tau(L,\rho)^2}{1-\rho^2}
\norm{A}_+^2\norm{P}_+^2\norm{B}_+^2\r^2$,
$C_3=3\frac{\tau(L,\rho)^2}{1-\rho^2}\norm{L}^2\norm{S}$, and
$C_4=44\frac{\tau(L,\rho)^2}{1-\rho^2}
\norm{P}_+^3\norm{A}_+^2\norm{B}_+^5\r^3$.
\end{lemma}

\begin{proof}
First we upper bound the operator norm of the linear operator $\T^\inv$.

\begin{lemma}
\label{lem:t_inv_operator_norm}
Let $\tau(M, \rho)$ denote the rate of growth of a matrix $M$. Then we have that
\(
\norm{\T^\inv}
\leq
\frac{\tau(L,\rho)^2}{1-\rho^2}.
\)
\end{lemma}

\begin{proof}
Since $L$ is a stable matrix, $\T$ is invertible. Moreover, whenever $L$ is stable and $X-L^\top XL=M$ for some matrix $M$, we know that
\(
X=\sum_{t\geq0}(L^t)^\top M L^t
\)
due to properties of solutions to Lyapunov equations~\citep{anderson2007optimal}. Hence we have that
\begin{align*}
\norm{\T^\inv}
&=
\sup_{\norm{M}=1}\norm{\T^\inv(M)} =
\sup_{\norm{M}=1}
\norm{\sum_{t\geq0}(L^t)^\top M L^t} \\
& \leq
\sum_{t\geq0}\norm{L^t}^2 \leq
\sum_{t\geq0}\tau(L,\rho)^2\rho^{2t}
=
\frac{\tau(L,\rho)^2}{1-\rho^2},
\end{align*}
where in the final line we used the definition of $\tau(L,\rho)$ and the geometric sum.
\end{proof}

Next we recall that
\[
\H(X)=L^\top X(I+S(P+X))^\inv SXL.
\]
Using \cref{lem:inverse} (appendix), we have that
\begin{align*}
\norm{\H(X)}
&=
\norm{L^\top X(I+S(P+X))^\inv SXL} \\
&\leq
\norm{L}^2\norm{S}\norm{X}^2.
\end{align*}

Let $P_X$ denote $P+X$ and consider the difference
\[
F(P_X,A,B,Q,R)
-
F(P_X,A_\epsilon,B_\epsilon,Q_\epsilon,R_\epsilon).
\]
Using \cref{eq:f_definition}, we have
\ifdefined\arxiv
\begin{align*}
F(P_X,A,B,Q,R)
-F(P_X,A_\epsilon,B_\epsilon,Q_\epsilon,R_\epsilon) =
A_\epsilon^\top P_X(I+S_\epsilon P_X)^{-1}A_\epsilon
-
A^\top P_X(I+SP_X)^{-1}A
+\Delta_Q .
\end{align*}
\else
\begin{align*}
&F(P_X,A,B,Q,R)
-F(P_X,A_\epsilon,B_\epsilon,Q_\epsilon,R_\epsilon) \\
&=
A_\epsilon^\top P_X(I+S_\epsilon P_X)^{-1}A_\epsilon
-
A^\top P_X(I+SP_X)^{-1}A
+\Delta_Q .
\end{align*}
\fi
Moreover, since $A_\epsilon=A+\Delta_A$, we have
\ifdefined\arxiv
\begin{align*}
F(P_X,A,B,Q,R)
-F(P_X,A_\epsilon,B_\epsilon,Q_\epsilon,R_\epsilon) & =
-A^\top P_X(I+SP_X)^{-1}
\Delta_S P_X(I+S_\epsilon P_X)^{-1}A \\
&\quad
+A^\top P_X(I+S_\epsilon P_X)^{-1}\Delta_A \\
&\quad
+\Delta_A^\top P_X(I+S_\epsilon P_X)^{-1}A
+\Delta_A^\top P_X(I+S_\epsilon P_X)^{-1}\Delta_A
+\Delta_Q .
\end{align*}

\else
\begin{align*}
&F(P_X,A,B,Q,R)
-F(P_X,A_\epsilon,B_\epsilon,Q_\epsilon,R_\epsilon) \\
&=
-A^\top P_X(I+SP_X)^{-1}
\Delta_S P_X(I+S_\epsilon P_X)^{-1}A \\
&\quad
+A^\top P_X(I+S_\epsilon P_X)^{-1}\Delta_A \\
&\quad
+\Delta_A^\top P_X(I+S_\epsilon P_X)^{-1}A \\
&\quad
+\Delta_A^\top P_X(I+S_\epsilon P_X)^{-1}\Delta_A
+\Delta_Q .
\end{align*}
\fi
By again using \cref{lem:inverse},
\ifdefined\arxiv
\begin{align*}
\norm{
F(P_X,A,B,Q,R)
-F(P_X,A_\epsilon,B_\epsilon,Q_\epsilon,R_\epsilon)
} \leq
\norm{A}^2\norm{P_X}^2\norm{\Delta_S}
+2\norm{A}\norm{P_X}\epsilon +\norm{P_X}\epsilon^2+\epsilon .
\end{align*}
\else
\begin{align*}
&\norm{
F(P_X,A,B,Q,R)
-F(P_X,A_\epsilon,B_\epsilon,Q_\epsilon,R_\epsilon)
} \\
&\leq
\norm{A}^2\norm{P_X}^2\norm{\Delta_S}
+2\norm{A}\norm{P_X}\epsilon +\norm{P_X}\epsilon^2+\epsilon .
\end{align*}
\fi
However, since $X\in\S_\nu$ we have $\norm{X}\leq\nu$ and $\norm{P_X}\leq\norm{P}+\nu$. Additionally, as $\nu\leq1$, we have $\norm{P_X}\leq1+\norm{P}$. Moreover, we also assume $\epsilon\leq1$. Therefore,
\ifdefined\arxiv
\begin{align*}
\norm{\Phi(X)}
&\leq
\frac{\tau(L,\rho)^2}{1-\rho^2}
\Big(
\norm{\H(X)} 
+
\norm{
F(P_X,A,B,Q,R)
-F(P_X,A_\epsilon,B_\epsilon,Q_\epsilon,R_\epsilon)
}
\Big) \\
&\leq
\frac{\tau(L,\rho)^2}{1-\rho^2}
\Big(
\norm{L}^2\norm{S}\nu^2 
+4\norm{A}^2\norm{P_X}^2
\norm{B}_+^2\r^2\epsilon 
+2\norm{A}\norm{P_X}\epsilon
+\norm{P_X}\epsilon^2+\epsilon
\Big) \\
&\leq
\frac{\tau(L,\rho)^2}{1-\rho^2}
\Big(
\norm{L}^2\norm{S}\nu^2
+8\norm{A}_+^2\norm{P}_+^2
\norm{B}_+^2\r^2\epsilon
\Big) \\
&=
C_1\nu^2+C_2\epsilon.
\end{align*}
\else
\begin{align*}
\norm{\Phi(X)}
&\leq
\frac{\tau(L,\rho)^2}{1-\rho^2}
\Big(
\norm{\H(X)} \\
&\quad
+
\norm{
F(P_X,A,B,Q,R)
-F(P_X,A_\epsilon,B_\epsilon,Q_\epsilon,R_\epsilon)
}
\Big) \\
&\leq
\frac{\tau(L,\rho)^2}{1-\rho^2}
\Big(
\norm{L}^2\norm{S}\nu^2 \\
&\quad
+4\norm{A}^2\norm{P_X}^2
\norm{B}_+^2\r^2\epsilon \\
&\quad
+2\norm{A}\norm{P_X}\epsilon
+\norm{P_X}\epsilon^2+\epsilon
\Big) \\
&\leq
\frac{\tau(L,\rho)^2}{1-\rho^2}
\Big(
\norm{L}^2\norm{S}\nu^2
+8\norm{A}_+^2\norm{P}_+^2
\norm{B}_+^2\r^2\epsilon
\Big) \\
&=
C_1\nu^2+C_2\epsilon.
\end{align*}
\fi

Next up we show the bound on $\norm{\Phi(X_1)-\Phi(X_2)}$. Denote
\[
\G(X)
=
F(P_X,A,B,Q,R)
-
F(P_X,A_\epsilon,B_\epsilon,Q_\epsilon,R_\epsilon).
\]
Note that
\ifdefined\arxiv
\begin{align*}
\norm{\Phi(X_1)-\Phi(X_2)}
&\leq
\norm{\T^\inv}
\Big(
\norm{\G(X_1)-\G(X_2)} 
+
\norm{\H(X_1)-\H(X_2)}
\Big).
\end{align*}
\else
\begin{align*}
\norm{\Phi(X_1)-\Phi(X_2)}
&\leq
\norm{\T^\inv}
\Big(
\norm{\G(X_1)-\G(X_2)} \\
&\quad
+
\norm{\H(X_1)-\H(X_2)}
\Big).
\end{align*}
\fi
However,
\(
\norm{\T^\inv}
\leq
\frac{\tau(L,\rho)^2}{1-\rho^2}
\)
from \cref{lem:t_inv_operator_norm}. We thus deal with the other two terms individually.
Set $M_i=(I+SP_{X_i})^{-1}$ for $i=1,2$. Since $P_{X_i}=P+X_i$ and if $\Delta_X=X_1-X_2$,
\[
M_1-M_2
=
M_1S(P_{X_2}-P_{X_1})M_2
=
M_1S(-\Delta_X)M_2.
\]
Thus,
\ifdefined\arxiv
\begin{align*}
\H(X_1)-\H(X_2)
&=
L^\top
\Big[
X_1M_1S(-\Delta_X)M_2SX_1 
+\Delta_X M_2 S X_1
+X_2M_2 S \Delta_X
\Big]L .
\end{align*}
\else
\begin{align*}
\H(X_1)-\H(X_2)
&=
L^\top
\Big[
X_1M_1S(-\Delta_X)M_2SX_1 \\
&\qquad
+\Delta_X M_2 S X_1
+X_2M_2 S \Delta_X
\Big]L .
\end{align*}
\fi
By \Cref{lem:inverse}, the identity $(I+YX)^{-1}Y=Y(I+XY)^{-1}$, $\nu\leq\norm{S}^{-1}$, we have $\norm{M_i S} \leq \norm{S}$ and,
\begin{align*}
\norm{\H(X_1)-\H(X_2)}
&\leq
\norm{L}^2
\left(\norm{S}^2\nu^2+2\norm{S}\nu\right)
\norm{\Delta_X} \\
&\leq
3\norm{L}^2\norm{S}\nu\norm{\Delta_X}.
\end{align*}

Lastly we need to deal with the $\G$ terms. We start off by noting that
\begin{align*}
\norm{(I+SP_X)^\inv}
&=
\norm{P_X^\inv P_X(I+SP_X)^\inv} \\
&\leq
\norm{P_X^\inv}\norm{P_X}
\leq
2\norm{P_X},
\end{align*}
where we use the fact that $\norm{X}\leq\nu\leq1/2$ and $P\succeq I$.

For compactness, define
\[
M_X:=P_X(I+SP_X)^{-1},
\,\,\,
M_X^\epsilon:=P_X(I+S_\epsilon P_X)^{-1}.
\]
Using the definition of $\G$, we can write
\ifdefined\arxiv
\begin{align*}
\G(X) =
-A^\top M_X\Delta_S M_X^\epsilon A
+A^\top M_X^\epsilon\Delta_A 
+\Delta_A^\top M_X^\epsilon A
+\Delta_A^\top M_X^\epsilon\Delta_A
+\Delta_Q .
\end{align*}
\else
\begin{align*}
\G(X)
&=
-A^\top M_X\Delta_S M_X^\epsilon A
+A^\top M_X^\epsilon\Delta_A \\
&\quad
+\Delta_A^\top M_X^\epsilon A
+\Delta_A^\top M_X^\epsilon\Delta_A
+\Delta_Q .
\end{align*}
\fi
Therefore,
\ifdefined\arxiv
\begin{align*}\G(X_1)-\G(X_2)
=
-A^\top g(\Delta_X,S,S_\epsilon)A 
+A^\top f(\Delta_X,S_\epsilon)\Delta_A 
+\Delta_A^\top f(\Delta_X,S_\epsilon)A
+\Delta_A^\top f(\Delta_X,S_\epsilon)\Delta_A,
\end{align*}
\else
\begin{align*}
& \G(X_1)-\G(X_2)
=
-A^\top g(\Delta_X,S,S_\epsilon)A \\
&\quad
+A^\top f(\Delta_X,S_\epsilon)\Delta_A \\
&\quad
+\Delta_A^\top f(\Delta_X,S_\epsilon)A
+\Delta_A^\top f(\Delta_X,S_\epsilon)\Delta_A,
\end{align*}
\fi
\noindent where
\begin{align*}
f(\Delta_X,S)
&:=
P_{X_1}(I+SP_{X_1})^{-1}
-
P_{X_2}(I+SP_{X_2})^{-1}, \\
g(\Delta_X,S,S_\epsilon)
&:=
M_{X_1}\Delta_S M_{X_1}^\epsilon
-
M_{X_2}\Delta_S M_{X_2}^\epsilon .
\end{align*}
Moreover,
\ifdefined\arxiv
\begin{align*}
f(\Delta_X,S)
&=
\Delta_X(I+SP_{X_1})^{-1} 
+P_{X_2}
\left[
(I+SP_{X_1})^{-1}
-
(I+SP_{X_2})^{-1}
\right] \\
&=
\Delta_X(I+SP_{X_1})^{-1} 
+P_{X_2}(I+SP_{X_2})^{-1}
S(-\Delta_X)(I+SP_{X_1})^{-1}.
\end{align*}
\else
\begin{align*}
f(\Delta_X,S)
&=
\Delta_X(I+SP_{X_1})^{-1} \\
&\quad
+P_{X_2}
\left[
(I+SP_{X_1})^{-1}
-
(I+SP_{X_2})^{-1}
\right] \\
&=
\Delta_X(I+SP_{X_1})^{-1} \\
&\quad
+P_{X_2}(I+SP_{X_2})^{-1}
S(-\Delta_X)(I+SP_{X_1})^{-1}.
\end{align*}
\fi
Thus,
\begin{align*}
\norm{f(\Delta_X,S)}
&\leq
2\norm{\Delta_X}\norm{P_{X_1}} +2\norm{\Delta_X}
\norm{P_{X_2}}\norm{S}\norm{P_{X_1}} \\
&\leq
4\norm{\Delta_X}\norm{P}_+^2\norm{S}_+ .
\end{align*}
Similarly,
\begin{align*}
g(\Delta_X,S,S_\epsilon)
&=
f(\Delta_X,S)\Delta_S M_{X_1}^\epsilon
+
M_{X_2}\Delta_S f(\Delta_X,S_\epsilon),
\end{align*}
and hence
\[
\norm{g(\Delta_X,S,S_\epsilon)}
\leq
\norm{\Delta_S}\norm{P}_+
\left(
\norm{f(\Delta_X,S)}
+
\norm{f(\Delta_X,S_\epsilon)}
\right).
\]
Therefore, after some cumbersome algebra we can show
\ifdefined\arxiv
\begin{align*}
\norm{\G(X_1)-\G(X_2)}
&\leq
\norm{A}^2\norm{g(\Delta_X,S,S_\epsilon)} 
+2\norm{A}\norm{\Delta_A}
\norm{f(\Delta_X,S_\epsilon)} 
+\norm{\Delta_A}^2
\norm{f(\Delta_X,S_\epsilon)} \\
&\leq
\norm{A}^2\norm{P}_+\norm{\Delta_S}
\Big(
\norm{f(\Delta_X,S)} 
+\norm{f(\Delta_X,S_\epsilon)}
\Big) 
+2\epsilon\norm{A}\norm{f(\Delta_X,S_\epsilon)} +\epsilon^2\norm{f(\Delta_X,S_\epsilon)} \\
&\leq
4\norm{A}^2\norm{P}_+^3
\norm{\Delta_X}\norm{\Delta_S}
\left(\norm{S}_+ + \norm{S_\epsilon}_+\right)
+8\epsilon\norm{A}\norm{P}_+^2
\norm{\Delta_X}\norm{S_\epsilon}_+ 
+4\epsilon^2\norm{\Delta_X}\norm{P}_+^2
\norm{S_\epsilon}_+ .
\end{align*}
\else
\begin{align*}
\norm{\G(X_1)-\G(X_2)}
&\leq
\norm{A}^2\norm{g(\Delta_X,S,S_\epsilon)} \\
&\quad
+2\norm{A}\norm{\Delta_A}
\norm{f(\Delta_X,S_\epsilon)} \\
&\quad
+\norm{\Delta_A}^2
\norm{f(\Delta_X,S_\epsilon)} \\
&\leq
\norm{A}^2\norm{P}_+\norm{\Delta_S}
\Big(
\norm{f(\Delta_X,S)} \\
&\quad
+\norm{f(\Delta_X,S_\epsilon)}
\Big) \\
&\quad
+2\epsilon\norm{A}\norm{f(\Delta_X,S_\epsilon)} +\epsilon^2\norm{f(\Delta_X,S_\epsilon)} \\
&\leq
4\norm{A}^2\norm{P}_+^3
\norm{\Delta_X}\norm{\Delta_S}
\left(\norm{S}_+ + \norm{S_\epsilon}_+\right) \\
&\quad
+8\epsilon\norm{A}\norm{P}_+^2
\norm{\Delta_X}\norm{S_\epsilon}_+ \\
&\quad
+4\epsilon^2\norm{\Delta_X}\norm{P}_+^2
\norm{S_\epsilon}_+ .
\end{align*}
\fi

Combining all of the different terms with $\epsilon \leq 1$, $\norm{S_\epsilon} \leq 1+\norm{S}$, and $\norm{\Delta_S}\leq4\epsilon\norm{B}_+^2\r^2$ yields
\[
\norm{\Phi(X_1)-\Phi(X_2)}
\leq
\norm{X_1-X_2}(C_3\nu+C_4\epsilon),
\]
where
$C_3 = 3\frac{\tau(L,\rho)^2}{1-\rho^2}\norm{L}^2\norm{S}
$
and
$ C_4 = 44\frac{\tau(L,\rho)^2}{1-\rho^2} \norm{P}_+^3\norm{A}_+^2\norm{B}_+^5\r^3.
$
\end{proof}

Using \cref{lem:contraction}, we are able to finish the proof of \cref{thm:dare_sensitivity} as follows.  We first establish that $\Phi$ maps $\S_\nu$ to $\S_\nu$.  To see this, note that by \cref{lem:contraction},
\[
\norm{\Phi(X)}\leq C_1\nu^2+C_2\epsilon.
\]
Setting $\nu=2C_2\epsilon$,
\[
\norm{\Phi(X)}
\leq
C_1\nu^2+\frac{1}{2}\nu
=
\nu\left(C_1\nu+\frac{1}{2}\right).
\]
However, $\epsilon$ is chosen such that $\epsilon\leq\frac{1}{4}C_1^\inv C_2^\inv$ so that $C_1\nu\leq\frac{1}{2}$.  Moreover, we also have that $\nu\leq\frac{1}{2}$ and $\nu\leq\norm{S}^\inv$ by the bound on $\epsilon$. Thus we find that $\norm{\Phi(X)}\leq\nu$ and so $\Phi$ maps $\S_\nu$ to itself.

Next we show that $\Phi$ is a contraction over $\S_\nu$. By the choice of $\epsilon$ and \cref{lem:contraction},
\begin{align*}
\norm{\Phi(X_1)-\Phi(X_2)}
&\leq
\norm{X_1-X_2}(C_3\nu+C_4\epsilon) \\
&\leq
2\max\{C_3\nu,C_4\epsilon\}
\norm{X_1-X_2} \\
&\leq
\frac{1}{4}\norm{X_1-X_2}.
\end{align*}
Hence, $\Phi$ has a fixed point in $\S_\nu$ by the Banach fixed point theorem since $\S_\nu$ is a closed set. However, as established before, the fixed point of $\Phi$ is precisely $\Delta_P$. Therefore, $\Delta_P$ is in $\S_\nu$ and hence $\norm{\Delta_P}\leq\nu=2C_2\epsilon$. The final bound follows by plugging in the appropriate constants and verifying that the upper bound on $\epsilon$ satisfies the earlier inequalities.
\end{proof}

\subsection{Upper Bounds on Control and Stability Margin}
\label{sec:dare_constants}

Unfortunately, \cref{thm:dare_sensitivity} does not immediately suffice for showing our algorithm uniformly approximates $\PF(\S)$.  Denote by $K$ the optimal control for $\loss_w(\cdot)$ for a given scalarization parameter $w \in \Delta([m])$ and by $K_\epsilon$ the optimal control for $\loss_{w_\epsilon}(\cdot)$ for $w_\epsilon \in N_\epsilon$ with $\norm{w - w_\epsilon}_1 \leq \epsilon$.  Let $P = \dare(A,B,Q_w,R_w)$ and $P_\epsilon = \dare(A,B,Q_{w_\epsilon}, R_{w_\epsilon})$.  A direct application of \cref{thm:dare_sensitivity} provides guarantees on $\norm{P - P_\epsilon}$ in terms of quantities such as $\norm{P}_+$, $\norm{A+BK}_+$, and $\tau(A+BK, \rho)$.  All of these parameters can depend on the choice of $w \in \Delta([m])$.  Thus, in order to obtain guarantees that are uniform over the scalarization simplex, we first introduce $w$-independent bounds on the Riccati solutions and controllers.  We then consider the stability margins across the Pareto front.

\begin{lemma}
\label{lem:upper_bounds}
Let $K_0$ be any stabilizing controller for $(A,B)$, let $L_0=A+BK_0$, and choose $\rho_0$ such that $\rho(L_0)<\rho_0<1$.  Then
\[
\Pmax
\leq
\frac{\tau(L_0,\rho_0)^2}{1-\rho_0^2}
\left(
\max_j\norm{Q_j}
+
\norm{K_0}^2\max_j\norm{R_j}
\right).
\]
Moreover,
\[
\Kmax
=
\sup_{w\in\Delta([m])}\norm{K_w}
\leq
m\max_j\norm{R_j^\inv}\norm{B}\norm{A}\Pmax.
\]
\end{lemma}
\noindent The proof is deferred to \cref{app:omitted_proofs}.  Taking, for instance, $K_0$ to be the single-objective LQR controller for the cost $(I_n,I_d)$ gives a computable bound in terms of the input matrices.

Next, we tackle the stability margin $\tau(A+BK, \rho)$ which appears in the bounds of \cref{thm:dare_sensitivity}.  When using it to establish the performance guarantee for our algorithm we will let $K = K_w$ for a particular choice of scalarization parameter $w \in \Delta([m])$.  Hence, we need to provide an upper bound on $\sup_{w \in \Delta([m])} \tau(A+BK_w, \rho)$.  First note that each scalarized LQR controller $K_w$ is stabilizing, so there exists a $\gamma$ such that $\rho(A+BK_w) < \gamma < 1$ for that fixed $w$ and hence $\tau(A+BK_w, \gamma) < \infty$.  One attempt might be to establish uniform stability bounds over {\em all} stabilizing controllers, i.e. $\sup_{K \in \S} \tau(A+BK, \gamma)$.  However, there may be a sequence of stabilizing controllers for which the closed-loop spectral radius converges to one.  We next establish using a continuity argument that it {\em is possible} to provide a uniform stability margin over the Pareto front since $\Delta([m])$ is closed and $\tau(A+BK_w, \rho)$ is continuous with respect to the scalarization parameter.

\begin{lemma}
\label{lem:uniform_control_pf}
There exists a $\gbar$ such that $\rho(A+BK_w) \leq \gbar < 1$ for all $K_w \in \PF(\S)$.  Defining
\begin{equation}
\label{eq:stability_margin}
\tbar = \sup_{w \in \Delta([m])} \, \tau(A+BK_w, \gbar),
\end{equation}
we also have that $\tbar < \infty$.
\end{lemma}

\begin{proof}
The map $w\mapsto \rho(A+BK_w)$ is continuous (see \cref{lem:eig_cont}), each $K_w$ is stabilizing by the single-objective LQR optimality conditions, and $\Delta([m])$ is compact. Hence
\[
\bar\rho
:=
\sup_{w\in\Delta([m])}\rho(A+BK_w)
<1.
\]
Choose any $\gbar$ such that $\bar\rho<\gbar<1$. Then $\rho(A+BK_w)<\gbar$ for every $w\in\Delta([m])$, and by Gelfand's formula $\tau(A+BK_w,\gbar)<\infty$ for every $w$.

It remains to show that the supremum is finite. For fixed $\gbar$, the map
\(
w\mapsto \tau(A+BK_w,\gbar)
\)
is locally bounded on $\Delta([m])$. Indeed, if $\rho(A+BK_w)<\gbar$, then by continuity of $K_w$ and of the spectral radius, the same strict inequality holds in a neighborhood of $w$. This gives a local uniform bound on
\(
\norm{(A+BK_{w'})^k}\gbar^{-k}.
\)
Since $\Delta([m])$ is compact, finitely many such neighborhoods cover it, and therefore $\tbar < \infty$.
\end{proof}

We now use \cref{lem:upper_bounds} and \cref{lem:uniform_control_pf} together with \cref{thm:dare_sensitivity} to provide an interpretable bound on the sensitivity of solutions to the discrete algebraic Riccati equation over the Pareto front which does not depend on the choice of linear scalarization parameter $w \in \Delta([m])$.

\begin{corollary}
\label{cor:dare_sensitivity_clean}
Consider arbitrary $w$ and $w_\epsilon$ in $\Delta([m])$ such that $\norm{w - w_\epsilon}_1 \leq \epsilon$.  Further denote $P = \dare(A,B,Q_w,R_w)$, $P_\epsilon = \dare(A_\epsilon,B_\epsilon, Q_{w_\epsilon}, R_{w_\epsilon})$, and $K=K_w$.
Suppose that \Cref{ass:stabilizable_observable} holds, $(A_\epsilon,B_\epsilon)$ is stabilizable, and $\max\{\norm{\Delta_A}, \norm{\Delta_B}\} \leq \epsilon$.  Denote
\ifdefined\arxiv
\[
\Gamma
=
\max\left\{
1+\Pmax,\,
1+\Kmax,\,
\norm{A}_+,\,
\norm{B}_+,
\max_j\norm{Q_j},\,
\max_j\norm{R_j},\,
1+m\max_j\norm{R_j^\inv}
\right\}.
\]
\else
\[
\Gamma
=
\max\left\{
\begin{array}{c}
1+\Pmax,\,
1+\Kmax,\,
\norm{A}_+,\,
\norm{B}_+, \\
\max_j\norm{Q_j},\,
\max_j\norm{R_j},\,
1+m\max_j\norm{R_j^\inv}
\end{array}
\right\}.
\] \fi
Then we have that
\[
\norm{P-P_\epsilon}
\leq
O(1)
\frac{\tbar^2}{1-\gbar^2}
\Gamma^8\epsilon
\]
so long as
\[
\epsilon
\leq
O(1)
\frac{(1-\gbar^2)^4}{\tbar^4}
\Gamma^{-17}.
\]
\end{corollary}

\begin{proof}
First note that
\begin{align*}
\norm{Q_w-Q_{w_\epsilon}}
&=
\norm{\sum_i (w_i-w_{\epsilon,i})Q_i} \\
&\leq
\norm{w-w_\epsilon}_1\max_j\norm{Q_j}
\leq
\epsilon\max_j\norm{Q_j}.
\end{align*}
Similarly, $\norm{R_w - R_{w_\epsilon}} \leq \epsilon \max_j \norm{R_j}$.  Both $(R_w, Q_w)$ and $(R_{w_\epsilon}, Q_{w_\epsilon})$ are also positive definite by \Cref{ass:stabilizable_observable}.  Furthermore, $\norm{P}_+ \leq \Gamma$ and
\[
\norm{L}
=
\norm{A+BK}
\leq
\norm{A}+\norm{B}\Kmax
\leq
\norm{A}_+\norm{B}_+\Kmax.
\]
Hence, using \cref{thm:dare_sensitivity}, we have that if
\ifdefined\arxiv
\begin{align*}
\epsilon
\max_j\max\{\norm{Q_j},\norm{R_j}\}
&\leq
O(1)
\frac{(1-\gbar^2)^4}{\tbar^4}\Gamma^{-16} \\
&\leq
O(1)
\frac{(1-\gbar^2)^4}{\tau(L,\gbar)^4}
\norm{A}_+^{-2}
\norm{P}_+^{-3} 
\norm{B}_+^{-4}
\r^{-3}
\norm{L}_+^{-2},
\end{align*} 
\else
\begin{align*}
\epsilon
\max_j\max\{\norm{Q_j},\norm{R_j}\}
&\leq
O(1)
\frac{(1-\gbar^2)^4}{\tbar^4}\Gamma^{-16} \\
&\leq
O(1)
\frac{(1-\gbar^2)^4}{\tau(L,\gbar)^4}
\norm{A}_+^{-2}
\norm{P}_+^{-3} \\
&\quad
\cdot
\norm{B}_+^{-4}
\r^{-3}
\norm{L}_+^{-2},
\end{align*}
\fi
then
\begin{align*}
\norm{P-P_\epsilon}
&\leq
O(1)
\frac{\tau(L,\gbar)^2}{1-\gbar^2}
\norm{A}_+^2\norm{P}_+^2
\norm{B}_+^2\r^2\epsilon \\
&\leq
O(1)
\frac{\tbar^2}{1-\gbar^2}
\Gamma^8\epsilon .
\end{align*}
The final bound follows from \cref{lem:upper_bounds}.
\end{proof}

%% file: r1_parts/main_proofs.tex
\section{Proofs of Main Results}
\label{sec:main_proofs}

\begin{proof}[Proof of Theorem~\ref{thm:pf_approx}]
Let $K \in \PF(\S)$ be arbitrary. By \cref{thm:sufficiency_linear}, there exists a $w \in \Delta([m])$ such that $K$ optimizes $\loss_w(\cdot)$. Since $N_\epsilon$ is an $\epsilon$-net of $\Delta([m])$, let $w_\epsilon \in N_\epsilon$ be such that $\norm{w-w_\epsilon}_1 \leq \epsilon$. Set $K_\epsilon$ as the control in $\PF_\epsilon(\S)$ which optimizes $\loss_{w_\epsilon}(\cdot)$. Denote
\[
P=\dare(A,B,Q_w,R_w),
\qquad
P_\epsilon=\dare(A,B,Q_{w_\epsilon},R_{w_\epsilon}).
\]
Applying \cref{cor:dare_sensitivity_clean} and the assumption on $\epsilon$,
\begin{equation}
\label{eq:dare_approx_proof}
\norm{P-P_\epsilon}
\leq
O(1)
\frac{\tbar^2}{1-\gbar^2}
\Gamma^8\epsilon .
\end{equation}
As shown in \cref{lem:control_close_and_tau} (appendix)
\begin{equation}
\label{eq:control_approx_proof}
\norm{K-K_\epsilon}
\leq
7\Gamma^4\norm{P-P_\epsilon}
\leq
O(1)
\frac{\tbar^2}{1-\gbar^2}
\Gamma^{12}\epsilon .
\end{equation}

We now prove the entry-wise error bound. Fix an objective $i\in[m]$. By \cref{lem:cost_difference},
\ifdefined\arxiv
\begin{align*}
\loss_i(K_\epsilon)-\loss_i(K) =
\sigma^2
\Big[
-2\Tr\left(
P_{K_\epsilon}^{L}
(K-K_\epsilon)^\top E_K
\right)
+
\Tr\left(
P_{K_\epsilon}^{L}
(K-K_\epsilon)^\top
(R_i+B^\top P_KB)
(K-K_\epsilon)
\right)
\Big],
\end{align*}
\else
\begin{align*}
&\loss_i(K_\epsilon)-\loss_i(K) \\
&=
\sigma^2
\Big[
-2\Tr\left(
P_{K_\epsilon}^{L}
(K-K_\epsilon)^\top E_K
\right) \\
&\qquad
+
\Tr\left(
P_{K_\epsilon}^{L}
(K-K_\epsilon)^\top
(R_i+B^\top P_KB)
(K-K_\epsilon)
\right)
\Big],
\end{align*}
\fi
where
\ifdefined\arxiv
\begin{align*}
P_{K_\epsilon}^{L} =
\dlyape(A+BK_\epsilon,I), 
\quad P_K =
\dlyape((A+BK)^\top,Q_i+K^\top R_iK),
\quad E_K =
(R_i+B^\top P_KB)K+B^\top P_KA .
\end{align*}
\else
\begin{align*}
P_{K_\epsilon}^{L}
& =
\dlyape(A+BK_\epsilon,I), \\
P_K
& =
\dlyape((A+BK)^\top,Q_i+K^\top R_iK), \\
E_K
& =
(R_i+B^\top P_KB)K+B^\top P_KA .
\end{align*}
\fi
Using the Von Neumann trace inequality, we obtain
\ifdefined\arxiv
\begin{align}
\label{eq:entrywise_loss_bound}
|\loss_i(K_\epsilon)-\loss_i(K)|
&\leq
2\sigma^2 n
\norm{P_{K_\epsilon}^{L}}
\norm{K-K_\epsilon}
\norm{E_K} 
+\sigma^2 n
\norm{K-K_\epsilon}^2
\norm{P_{K_\epsilon}^{L}}
\norm{R_i+B^\top P_KB}.
\end{align}
\else
\begin{align}
\label{eq:entrywise_loss_bound}
|\loss_i(K_\epsilon)-\loss_i(K)|
&\leq
2\sigma^2 n
\norm{P_{K_\epsilon}^{L}}
\norm{K-K_\epsilon}
\norm{E_K} \\
&\,\,
+\sigma^2 n
\norm{K-K_\epsilon}^2
\norm{P_{K_\epsilon}^{L}}
\norm{R_i+B^\top P_KB}. \notag
\end{align}
\fi
From \cref{eq:control_approx_proof},
\(
\norm{K-K_\epsilon}
\leq
O(1)
\frac{\tbar^2}{1-\gbar^2}
\Gamma^{12}\epsilon .
\)
Moreover, by \cref{lem:lyapunov_bound},
\[
\norm{P_{K_\epsilon}^{L}}
\leq
\frac{\tau(A+BK_\epsilon,\gbar)^2}{1-\gbar^2}
\leq
\frac{\tbar^2}{1-\gbar^2},
\]
since $K_\epsilon \in \PF(\S)$. Next, we note that
\(
\norm{E_K}
\leq
2\Gamma^3\norm{P_K}.
\)
Again using \cref{lem:lyapunov_bound}, we have
\(
\norm{P_K}
\leq
\Gamma^3
\frac{\tbar^2}{1-\gbar^2}.
\)
Thus,
\(
\norm{E_K}
\leq
O(1)
\Gamma^6
\frac{\tbar^2}{1-\gbar^2}.
\)
Similarly,
\[
\norm{R_i+B^\top P_KB}
\leq
O(1)
\Gamma^5
\frac{\tbar^2}{1-\gbar^2}.
\]
Plugging these bounds into \cref{eq:entrywise_loss_bound}, we obtain
\ifdefined\arxiv
\begin{align*}
|\loss_i(K_\epsilon)-\loss_i(K)|
\leq
O(1)
\sigma^2 n
\Gamma^{18}
\frac{\tbar^6}{(1-\gbar^2)^3}
\epsilon 
+
O(1)
\sigma^2 n
\Gamma^{29}
\frac{\tbar^8}{(1-\gbar^2)^4}
\epsilon^2 .
\end{align*}
\else
\begin{align*}
|\loss_i(K_\epsilon)-\loss_i(K)|
&\leq
O(1)
\sigma^2 n
\Gamma^{18}
\frac{\tbar^6}{(1-\gbar^2)^3}
\epsilon \\
&\quad
+
O(1)
\sigma^2 n
\Gamma^{29}
\frac{\tbar^8}{(1-\gbar^2)^4}
\epsilon^2 .
\end{align*}
\fi
Since $\epsilon$ satisfies the assumed upper bound, the second term is dominated by the stated bound. Therefore,
\[
|\loss_i(K_\epsilon)-\loss_i(K)|
\leq
O(1)
\sigma^2 n
\Gamma^{29}
\frac{\tbar^8}{(1-\gbar^2)^4}
\epsilon .
\]
Since this holds for every $i\in[m]$, we conclude that
\[
\norm{\vec{\loss}(K)-\vec{\loss}(K_\epsilon)}_\infty
\leq
O(1)
\sigma^2 n
\Gamma^{29}
\frac{\tbar^8}{(1-\gbar^2)^4}
\epsilon .
\]
\end{proof}

\begin{proof}[Proof of \cref{thm:ce_pf_approximation}]
We start with the following lemma establishing that replacing $(A,B)$ with $(\Ahat,\Bhat)$ yields $O(\epsilon)$ approximations for the optimal controllers with fixed cost matrices.

\begin{lemma}
\label{lem:ce_approx}
Under the hypotheses of \cref{thm:ce_pf_approximation}, for all $w\in\Delta([m])$, suppose that $K$ optimizes $\loss_w(\cdot)$ and $\Khat$ optimizes $\hatloss_w(\cdot)$. Then $\Khat$ stabilizes $(A,B)$ and for every $i\in[m]$,
\[
|\loss_i(K)-\loss_i(\Khat)|
\leq
O(1)
\sigma^2 n
\Gamma^{29}
\frac{\tbar^8}{(1-\gbar^2)^4}
\epsilon .
\]
\end{lemma}

\begin{proof}
Using \cref{cor:dare_sensitivity_clean} to establish sensitivity of the solutions to
$\dare(A,B,Q_w,R_w)$ and 
$\dare(\Ahat,\Bhat,Q_w,R_w),$
as well as \cref{lem:control_close_and_tau}, we have
\begin{align}
\label{eq:ce_control_close}
\norm{K-\Khat}
&\leq
O(1)
\Gamma^{12}
\frac{\tbar^2}{1-\gbar^2}
\epsilon, \\
\label{eq:ce_tau_bound}
\tau(A+B\Khat,(1+\gbar)/2)
&\leq
O(1)\tau(A+BK,\gbar).
\end{align}
In particular, $\Khat$ stabilizes $(A,B)$.
Using \cref{lem:cost_difference}, we have
\ifdefined\arxiv
\begin{align*}
\loss_i(\Khat)-\loss_i(K) =
\sigma^2
\Big[
-2\Tr\left(
P_{\Khat}^{L}
(K-\Khat)^\top E_K
\right) 
+
\Tr\left(
P_{\Khat}^{L}
(K-\Khat)^\top
(R_i+B^\top P_KB)
(K-\Khat)
\right)
\Big],
\end{align*}
\else
\begin{align*}
&\loss_i(\Khat)-\loss_i(K) \\
&=
\sigma^2
\Big[
-2\Tr\left(
P_{\Khat}^{L}
(K-\Khat)^\top E_K
\right) \\
&\qquad
+
\Tr\left(
P_{\Khat}^{L}
(K-\Khat)^\top
(R_i+B^\top P_KB)
(K-\Khat)
\right)
\Big],
\end{align*}
\fi
where
\ifdefined\arxiv
\begin{align*}
P_{\Khat}^{L}  =\dlyape(A+B\Khat,I), \quad
P_K  = \dlyape((A+BK)^\top,Q_i+K^\top R_iK), \quad
E_K  =(R_i+B^\top P_KB)K+B^\top P_KA.
\end{align*}
\else
\begin{align*}
P_{\Khat}^{L} & =\dlyape(A+B\Khat,I), \\
P_K & = \dlyape((A+BK)^\top,Q_i+K^\top R_iK), \\
E_K & =(R_i+B^\top P_KB)K+B^\top P_KA.
\end{align*}
\fi
Hence, by the Von Neumann trace inequality, denoting $\Delta_K=\Khat-K$,
\ifdefined\arxiv
\begin{align}
\label{eq:ce_entrywise_bound}
|\loss_i(\Khat)-\loss_i(K)|
&\leq
2\sigma^2 n
\norm{P_{\Khat}^{L}}
\norm{\Delta_K}
\norm{E_K} 
+\sigma^2 n
\norm{\Delta_K}^2
\norm{P_{\Khat}^{L}}
\norm{R_i+B^\top P_KB}. 
\end{align}
\else
\begin{align}
\label{eq:ce_entrywise_bound}
|\loss_i(\Khat)-\loss_i(K)|
&\leq
2\sigma^2 n
\norm{P_{\Khat}^{L}}
\norm{\Delta_K}
\norm{E_K} \\
&\quad
+\sigma^2 n
\norm{\Delta_K}^2
\norm{P_{\Khat}^{L}}
\norm{R_i+B^\top P_KB}. \notag
\end{align}
\fi
We deal with the terms one by one. First, notice that
\(
\norm{E_K}
\leq
2\Gamma^3\norm{P_K}.
\)
By \cref{lem:lyapunov_bound},
\(
\norm{P_K}
\leq
\Gamma^3
\frac{\tbar^2}{1-\gbar^2}.
\)
Thus,
\(
\norm{E_K}
\leq
O(1)
\Gamma^6
\frac{\tbar^2}{1-\gbar^2}.
\)
Moreover, by \cref{eq:ce_tau_bound} and \cref{lem:lyapunov_bound},
\(
\norm{P_{\Khat}^{L}}
\leq
O(1)
\frac{\tbar^2}{1-\gbar^2}.
\)
Similarly,
\[
\norm{R_i+B^\top P_KB}
\leq
O(1)
\Gamma^5
\frac{\tbar^2}{1-\gbar^2}.
\]
Combining these bounds with \cref{eq:ce_control_close} and \cref{eq:ce_entrywise_bound}, we obtain
\ifdefined\arxiv
\begin{align*}
|\loss_i(\Khat)-\loss_i(K)|
&\leq
O(1)
\sigma^2 n
\Gamma^{18}
\frac{\tbar^6}{(1-\gbar^2)^3}
\epsilon 
+
O(1)
\sigma^2 n
\Gamma^{29}
\frac{\tbar^8}{(1-\gbar^2)^4}
\epsilon^2 .
\end{align*}
\else
\begin{align*}
|\loss_i(\Khat)-\loss_i(K)|
&\leq
O(1)
\sigma^2 n
\Gamma^{18}
\frac{\tbar^6}{(1-\gbar^2)^3}
\epsilon \\
&\quad
+
O(1)
\sigma^2 n
\Gamma^{29}
\frac{\tbar^8}{(1-\gbar^2)^4}
\epsilon^2 .
\end{align*}
\fi
Since $\epsilon$ satisfies the assumed upper bound, the second term is dominated by the stated bound. Therefore,
\[
|\loss_i(\Khat)-\loss_i(K)|
\leq
O(1)
\sigma^2 n
\Gamma^{29}
\frac{\tbar^8}{(1-\gbar^2)^4}
\epsilon .
\]
\end{proof}

With the previous lemma in hand, \cref{thm:ce_pf_approximation} follows via the triangle inequality, \cref{lem:ce_approx}, and \cref{thm:pf_approx}. Let $K\in\PF(\S)$ be arbitrary and let $w\in\Delta([m])$ be such that $K$ optimizes $\loss_w(\cdot)$. Since $N_\epsilon$ is an $\epsilon$-net of $\Delta([m])$, let $w_\epsilon\in N_\epsilon$ satisfy $\norm{w-w_\epsilon}_1\leq\epsilon$. Denote by $\Khat_\epsilon$ the control which optimizes $\hatloss_{w_\epsilon}(\cdot)$ and by $K_\epsilon$ the control which optimizes $\loss_{w_\epsilon}(\cdot)$. Note that $K_\epsilon\in\PF_\epsilon(\S)$ and $\Khat_\epsilon\in\hatPF_\epsilon(\S)$. Also, by \cref{lem:ce_approx}, $\Khat_\epsilon$ stabilizes $(A,B)$.

For each $i\in[m]$,
\ifdefined\arxiv
\begin{align*}
|\loss_i(K)-\loss_i(\Khat_\epsilon)|
&\leq
|\loss_i(K)-\loss_i(K_\epsilon)| 
+
|\loss_i(K_\epsilon)-\loss_i(\Khat_\epsilon)|.
\end{align*}
\else
\begin{align*}
|\loss_i(K)-\loss_i(\Khat_\epsilon)|
&\leq
|\loss_i(K)-\loss_i(K_\epsilon)| \\
&\quad
+
|\loss_i(K_\epsilon)-\loss_i(\Khat_\epsilon)|.
\end{align*}
\fi
The first term is bounded by \cref{thm:pf_approx}, and the second term is bounded by \cref{lem:ce_approx}. Therefore,
\[
|\loss_i(K)-\loss_i(\Khat_\epsilon)|
\leq
O(1)
\sigma^2 n
\Gamma^{29}
\frac{\tbar^8}{(1-\gbar^2)^4}
\epsilon .
\]
Since this holds for all $i\in[m]$, we conclude that
\[
\norm{\vec{\loss}(K)-\vec{\loss}(\Khat_\epsilon)}_\infty
\leq
O(1)
\sigma^2 n
\Gamma^{29}
\frac{\tbar^8}{(1-\gbar^2)^4}
\epsilon .
\]
\end{proof}

%% file: r1_parts/conclusion.tex
\section{Conclusion}

We studied Pareto-front approximation for multi-objective LQR. Building on the classical scalarization characterization, we showed that a finite grid of scalarization weights yields a uniform, objective-wise approximation of the Pareto front. We also extended the guarantee to certainty-equivalent control with unknown dynamics.  Several directions remain open. First, the constants and polynomial dependence in our approximation bounds are unlikely to be tight, both in the known-dynamics and certainty-equivalent settings.  Second, the algorithm requires a number of LQR solves that is exponential in the number of objectives. Developing adaptive discretization or sampling methods that exploit the geometry of the Pareto front could reduce this dependence. 

%% file: r1_parts/appendix_lemmas.tex
\section{Auxiliary Lemmas and Omitted Proofs}
\label{app:omitted_proofs}

In this section we collect the omitted proof of \cref{lem:upper_bounds} and the auxiliary lemmas used throughout the proofs of \cref{thm:dare_sensitivity,thm:pf_approx,thm:ce_pf_approximation}.  We begin with the proof of \cref{lem:upper_bounds}.

\begin{proof}
We first prove \cref{lem:upper_bounds}. Fix a stabilizing controller $K_0$ and set $L_0=A+BK_0$.  For any $w\in\Delta([m])$, let $P_w=\dare(A,B,Q_w,R_w)$ and let
\[
P_{0,w}
=
\dlyape(L_0^\top,Q_w+K_0^\top R_wK_0)
\]
be the value matrix obtained by using the controller $K_0$ on the scalarized LQR instance.  Since $K_w$ is optimal for the scalarized problem, the value under $K_w$ is no larger than the value under $K_0$.  Hence $P_w\preceq P_{0,w}$ and
\[
\norm{P_w}\leq \norm{P_{0,w}}.
\]
By \cref{lem:lyapunov_bound},
\[
\norm{P_{0,w}}
\leq
\frac{\tau(L_0,\rho_0)^2}{1-\rho_0^2}
\norm{Q_w+K_0^\top R_wK_0}.
\]
Using $w\in\Delta([m])$,
\[
\norm{Q_w+K_0^\top R_wK_0}
\leq
\max_j\norm{Q_j}
+
\norm{K_0}^2\max_j\norm{R_j}.
\]
Taking the supremum over $w\in\Delta([m])$ gives the stated bound on $\Pmax$.

Next consider
\[
\Kmax
=
\sup_{w\in\Delta([m])}\norm{K_w}.
\]
For any $w\in\Delta([m])$, the controller $K_w$ optimizes $\loss_w(\cdot)$ and
\[
K_w
=
-(R_w+B^\top P B)^\inv B^\top P A,
\qquad
P=\dare(A,B,Q_w,R_w).
\]
Hence,
\begin{align*}
\sup_{w\in\Delta([m])}\norm{K_w}
&\leq
\sup_{w\in\Delta([m])}
\norm{R_w^\inv}\norm{B}\norm{A}\norm{P} \\
&\leq
m\max_j\norm{R_j^\inv}
\norm{B}\norm{A}\Pmax,
\end{align*}
where the first inequality uses \cref{lem:upper_bound_norm_sum} and the final inequality follows from the finiteness of $\Pmax$ and \cref{lem:bound_weighted_sum_eigenvalues}. Thus $\Kmax<\infty$.
\end{proof}

\subsection{Matrix Properties}
\label{app:matrix_lemmas}

\begin{lemma}[Lemma 7 of \cite{mania2019certainty}]
\label{lem:inverse}
Let $M$ and $N$ be two positive semidefinite matrices of the same dimension. Then
\[
\norm{N(I+MN)^\inv}\leq\norm{N}.
\]
\end{lemma}

\begin{lemma}
\label{lem:upper_bound_norm_sum}
Suppose that $A$ is positive definite, $B$ is positive semidefinite, and $A+B$ is invertible. Then
\[
\norm{(A+B)^\inv}\leq\norm{A^\inv}.
\]
\end{lemma}

\begin{proof}
Since $A$ and $B$ are positive semidefinite, $\smin(A)\leq\smin(A+B)$. Hence,
\[
\frac{1}{\smin(A+B)}
\leq
\frac{1}{\smin(A)},
\]
which is equivalent to
\[
\smax((A+B)^\inv)\leq\smax(A^\inv).
\]
Thus, $\norm{(A+B)^\inv}\leq\norm{A^\inv}$.
\end{proof}

\begin{lemma}
\label{lem:bound_weighted_sum_eigenvalues}
Let $w\in\Delta([m])$ and suppose that $R_w=\sum_i w_iR_i$, where each $R_i$ is positive definite. Then
\[
\norm{R_w^\inv}
\leq
m\max_j\norm{R_j^\inv}.
\]
\end{lemma}

\begin{proof}
For every $j\in[m]$ with $w_j \neq 0$,
\(
w_jR_j\preceq R_w.
\)
Thus, by \cref{lem:upper_bound_norm_sum},
\[
\norm{R_w^\inv}
\leq
\norm{(w_jR_j)^\inv}
=
\frac{1}{w_j}\norm{R_j^\inv}.
\]
Since $w\in\Delta([m])$, there exists an index $j$ with $w_j\geq 1/m$. For this index,
\(
\norm{R_w^\inv}
\leq
m\max_j\norm{R_j^\inv}.
\)
\end{proof}

\begin{lemma}
\label{lem:small_eigenvalue_bound}
Let $w\in\Delta([m])$ and suppose that $R_w=\sum_i w_iR_i$, where each $R_i$ is positive definite. Then
\ifdefined\arxiv
\begin{align*}
\lambda_n\bigl(B\bigl(\sum_i R_i\bigr)^\inv B^\top\bigr) 
 \leq
\lambda_n(BR_w^\inv B^\top)
 \leq
m\max_j\lambda_n(BR_j^\inv B^\top).
\end{align*}
\else
\begin{align*}
\lambda_n\bigl(B\bigl(\sum_i R_i\bigr)^\inv B^\top\bigr) 
& \leq
\lambda_n(BR_w^\inv B^\top) \\
& \leq
m\max_j\lambda_n(BR_j^\inv B^\top).
\end{align*}
\fi
\end{lemma}

\begin{proof}
For every $j\in[m]$,
\(
w_jR_j\preceq R_w\preceq \sum_i R_i.
\)
Therefore,
\[
B\Big(\sum_i R_i\Big)^\inv B^\top
\preceq
BR_w^\inv B^\top
\preceq
\frac{1}{w_j}BR_j^\inv B^\top.
\]
Since $w\in\Delta([m])$, there exists an index $j$ with $w_j\geq 1/m$. This gives the desired bound.
\end{proof}

\subsection{LQR Properties}
\label{app:lqr_lemmas}

\begin{lemma}
\label{lem:eigenvalue_dare_bound}
Let $w\in\Delta([m])$ and set $P_w=\dare(A,B,Q_w,R_w)$. If $\smin(Q_i)\geq 1$ for all $i\in[m]$, then
$\smin(P_w)\geq 1.$
\end{lemma}

\begin{proof}
This follows directly from the lower bounds on $\dare(\cdot)$ presented in \cite{garloff1986bounds}.
\end{proof}

\begin{lemma}
\label{lem:eig_cont}
Consider the map $f:\Delta([m])\to\mathbb{R}_{\geq0}$ defined by
\(
f(w)=\rho(A+BK_w),
\)
where $K_w$ optimizes $\loss_w(\cdot)$. Then $f$ is continuous.
\end{lemma}

\begin{proof}
First note that
\[
K_w
=
-(R_w+B^\top P_wB)^\inv B^\top P_wA,
\]
where $P_w=\dare(A,B,Q_w,R_w).$  Thus,
\[
f(w)
=
\rho\left(
A-B(R_w+B^\top P_wB)^\inv B^\top P_wA
\right).
\]
The map $M\mapsto\rho(M)$ is continuous. Moreover, $P_w=\dare(A,B,Q_w,R_w)$ is continuous with respect to $w$ by Lemma 2.1 in \cite{stoorvogel1996continuity}, using that $\smin(R_w)>0$ for all $w\in\Delta([m])$. Hence, $f$ is continuous by composition.
\end{proof}

\begin{lemma}
\label{lem:lyapunov_bound}
Suppose that $P=\dlyape(A,Q)$ with $A$ Schur stable. Then
\[
\norm{P}
\leq
\frac{\tau(A,\rho)^2}{1-\rho^2}\norm{Q}.
\]
\end{lemma}

\begin{proof}
Since $A$ is Schur stable,
\(
P=\sum_{t\geq0}(A^\top)^tQA^t.
\)
Therefore,
\begin{align*}
\norm{P}
&\leq
\norm{Q}\sum_{t\geq0}\norm{A^t}^2 \\
&\leq
\norm{Q}\sum_{t\geq0}\tau(A,\rho)^2\rho^{2t}
=
\norm{Q}\frac{\tau(A,\rho)^2}{1-\rho^2}.
\end{align*}
\end{proof}

\begin{lemma}
\label{lem:three_objectives}
For any stable $K\in\S$,
\[
\loss(K,Q,R)
=
\sigma^2\Tr(P_K)
=
\sigma^2\Tr\big((Q+K^\top RK)P_K^L\big),
\]
where $P_K=\dlyape((A+BK)^\top,Q+K^\top RK),$ and $P_K^L=\dlyape(A+BK,I).$
Moreover, if $K$ minimizes $\loss(\cdot,Q,R)$, then $ \loss(K,Q,R)=\sigma^2\Tr(P)$ where
$P=\dare(A,B,Q,R).$
\end{lemma}

\begin{proof}
Let $A_K = A + BK$ and $C_K = Q + K^\top R K$. Since $K$ is stable, the Lyapunov equation
\[
P_K = C_K + A_K^\top P_K A_K
\]
has a unique positive semidefinite solution, namely $P_K = \dlyape(A_K^\top,C_K)$. Therefore,
\[
x_t^\top C_K x_t
=
x_t^\top P_K x_t
-
x_t^\top A_K^\top P_K A_K x_t.
\]
Using $x_{t+1} = A_K x_t + \xi_t$ and $\Exp{\xi_t \xi_t^\top} = \sigma^2 I$, we also have
\[
\Exp{x_{t+1}^\top P_K x_{t+1} \mid x_t}
=
x_t^\top A_K^\top P_K A_K x_t
+
\sigma^2 \Tr(P_K).
\]
Combining the previous two displays gives
\[
\Exp{x_t^\top C_K x_t}
=
\Exp{x_t^\top P_K x_t}
-
\Exp{x_{t+1}^\top P_K x_{t+1}}
+
\sigma^2 \Tr(P_K).
\]
Averaging over $t = 0,\ldots,T$ and taking $T \to \infty$, the telescoping term vanishes because $A_K$ is stable and the state covariance remains uniformly bounded. Hence,
\[
\loss(K,Q,R)
=
\sigma^2 \Tr(P_K).
\]

Next, using the series representation of the Lyapunov equation and the cyclic property of trace,
\begin{align*}
\loss(K,Q,R)
&=
\sigma^2 \Tr(P_K) \\
&=
\sigma^2 \Tr\biggl(
\sum_{k \geq 0}
(A_K^\top)^k C_K A_K^k
\biggr) \\
&=
\sigma^2 \Tr\biggl(
C_K
\sum_{k \geq 0}
A_K^k (A_K^\top)^k
\biggr) \\
&=
\sigma^2 \Tr(C_K P_K^L).
\end{align*}
Since $C_K = Q + K^\top R K$ and $P_K^L = \dlyape(A_K,I)$, this proves the second equality.

Finally, if $K$ minimizes $\loss(\cdot,Q,R)$, then the standard LQR optimality equations imply that its value matrix $P$ solves
\[
P
=
Q + A^\top P A
-
A^\top P B(R + B^\top P B)^\inv B^\top P A.
\]
Thus $P = \dare(A,B,Q,R)$. Applying the first part of the proof to the optimal controller gives
\[
\loss(K,Q,R)
=
\sigma^2 \Tr(P).
\]
\end{proof}

\begin{lemma}
\label{lem:dare_equivalence}
Let $R$ be positive definite and let $R=LL^\top$ be its Cholesky decomposition. Then
\[
\dare(A,B,Q,R)
=
\dare(A,BL^{-\top},Q,I).
\]
\end{lemma}

\begin{proof}
Let $P=\dare(A,B,Q,R)$ and set $\tilde B=BL^{-\top}$. Since
\[
P
=
A^\top PA
-
A^\top PB(R+B^\top PB)^\inv B^\top PA
+
Q,
\]
it suffices to show that
\[
B(R+B^\top PB)^\inv B^\top
=
\tilde B(I+\tilde B^\top P\tilde B)^\inv\tilde B^\top.
\]
Indeed,
\begin{align*}
B(R+B^\top PB)^\inv B^\top
&=
B(R(I+R^\inv B^\top PB))^\inv B^\top \\
&=
B(I+R^\inv B^\top PB)^\inv R^\inv B^\top \\
&=
BL^{-\top}
(I+L^{-1}B^\top PBL^{-\top})^\inv
L^{-1}B^\top \\
&=
\tilde B(I+\tilde B^\top P\tilde B)^\inv\tilde B^\top,
\end{align*}
where the last step uses $(I+AB)^\inv A=A(I+BA)^\inv$.
\end{proof}

\begin{lemma}
\label{lem:cost_difference}
Fix $i\in[m]$ and let $K,K'\in\S$ be stable controllers.
Then
\ifdefined\arxiv
\begin{align*}
\loss_i(K')-\loss_i(K) =
\sigma^2
\Big[
2\Tr\left(
P_{K'}^L(K'-K)^\top E_K
\right) 
+
\Tr\left(
P_{K'}^L(K'-K)^\top
(R_i+B^\top P_KB)
(K'-K)
\right)
\Big].
\end{align*}

\else
\begin{align*}
&\loss_i(K')-\loss_i(K) \\
&=
\sigma^2
\Big[
2\Tr\left(
P_{K'}^L(K'-K)^\top E_K
\right) \\
&\quad
+
\Tr\left(
P_{K'}^L(K'-K)^\top
(R_i+B^\top P_KB)
(K'-K)
\right)
\Big].
\end{align*}
\fi
where $P_K = \dlyape((A+BK)^\top,Q_i+K^\top R_iK)$, $P_{K'}^L = \dlyape(A+BK',I),$ and $E_K = (R_i+B^\top P_KB)K+B^\top P_KA.$
\end{lemma}

\begin{proof}
Let $A_K = A+BK$, $A_{K'} = A+BK'$, and $C_K = Q_i + K^\top R_i K$.
Since $K$ is stable, $P_K = C_K + A_K^\top P_K A_K$.
Thus by \cref{lem:three_objectives},
\[
\loss_i(K)=\sigma^2\Tr(P_K),
\qquad
\loss_i(K')=\sigma^2\Tr(P_{K'}),
\]
where $P_{K'} = \dlyape(A_{K'}^\top,Q_i+K'^\top R_iK').$
Thus,
\[
\loss_i(K')-\loss_i(K)
=
\sigma^2\Tr(P_{K'}-P_K).
\]

Let $\Delta_K=K'-K$. We first derive a Lyapunov equation for
$P_{K'}-P_K$. Since
\[
P_{K'}
=
Q_i+K'^\top R_iK'
+
A_{K'}^\top P_{K'}A_{K'},
\]
we have
\ifdefined\arxiv
\begin{align*}
P_{K'}-P_K
=
A_{K'}^\top(P_{K'}-P_K)A_{K'} 
+
\Big[
Q_i+K'^\top R_iK'
+
A_{K'}^\top P_KA_{K'} 
-
Q_i-K^\top R_iK
-
A_K^\top P_KA_K
\Big].
\end{align*}
\else
\begin{align*}
P_{K'}-P_K
&=
A_{K'}^\top(P_{K'}-P_K)A_{K'} \\
&\quad
+
\Big[
Q_i+K'^\top R_iK'
+
A_{K'}^\top P_KA_{K'} \\
&\qquad
-
Q_i-K^\top R_iK
-
A_K^\top P_KA_K
\Big].
\end{align*}
\fi
Define the bracketed term by $H$. Then
\[
P_{K'}-P_K
=
\dlyape(A_{K'}^\top,H).
\]
Therefore, using the series representation of the Lyapunov equation and
cyclicity of trace,
\[
\Tr(P_{K'}-P_K)
=
\Tr(P_{K'}^L H),
\]
where $P_{K'}^L = \dlyape(A_{K'}, I)$.

It remains to expand $H$. Since $K'=K+\Delta_K$,
\[
K'^\top R_iK'-K^\top R_iK
=
\Delta_K^\top R_i\Delta_K
+
\Delta_K^\top R_iK
+
K^\top R_i\Delta_K .
\]
Also, since $A_{K'}=A_K+B\Delta_K$,
\ifdefined\arxiv
\begin{align*}
A_{K'}^\top P_KA_{K'}-A_K^\top P_KA_K
=
\Delta_K^\top B^\top P_KB\Delta_K 
+\Delta_K^\top B^\top P_KA_K
+
A_K^\top P_KB\Delta_K .
\end{align*}
\else
\begin{align*}
A_{K'}^\top P_KA_{K'}-A_K^\top P_KA_K
&=
\Delta_K^\top B^\top P_KB\Delta_K \\
&\quad
+\Delta_K^\top B^\top P_KA_K
+
A_K^\top P_KB\Delta_K .
\end{align*}
\fi
Combining the two displays gives
\[
H
=
\Delta_K^\top(R_i+B^\top P_KB)\Delta_K
+
\Delta_K^\top E_K
+
E_K^\top\Delta_K.
\]
Hence,
\ifdefined\arxiv
\begin{align*}
\Tr(P_{K'}^L H)
=
\Tr\left(
P_{K'}^L\Delta_K^\top(R_i+B^\top P_KB)\Delta_K
\right) 
+
\Tr\left(P_{K'}^L\Delta_K^\top E_K\right)
+
\Tr\left(P_{K'}^LE_K^\top\Delta_K\right).
\end{align*}
\else
\begin{align*}
\Tr(P_{K'}^L H)
&=
\Tr\left(
P_{K'}^L\Delta_K^\top(R_i+B^\top P_KB)\Delta_K
\right) \\
&\quad
+
\Tr\left(P_{K'}^L\Delta_K^\top E_K\right)
+
\Tr\left(P_{K'}^LE_K^\top\Delta_K\right).
\end{align*}
\fi
Since $P_{K'}^L$ is symmetric, the last two trace terms are equal. Thus,
\ifdefined\arxiv
\begin{align*}
\Tr(P_{K'}^L H)
=
2\Tr\left(P_{K'}^L\Delta_K^\top E_K\right) 
+
\Tr\left(
P_{K'}^L\Delta_K^\top(R_i+B^\top P_KB)\Delta_K
\right).
\end{align*}
\else
\begin{align*}
\Tr(P_{K'}^L H)
&=
2\Tr\left(P_{K'}^L\Delta_K^\top E_K\right) \\
&\quad
+
\Tr\left(
P_{K'}^L\Delta_K^\top(R_i+B^\top P_KB)\Delta_K
\right).
\end{align*}
\fi
Multiplying by $\sigma^2$ gives the bound.
\end{proof}


\begin{lemma}[Lemma 5 of \cite{mania2019certainty}]
\label{lem:tau_bound_sum}
Let $M\in\RR^{n\times n}$ and let $\rho\geq\rho(M)$. Then, for all $k\geq1$ and all matrices $\Delta$ of appropriate dimensions,
\[
\norm{(M+\Delta)^k}
\leq
\tau(M,\rho)
\left(
\tau(M,\rho)\norm{\Delta}+\rho
\right)^k.
\]
\end{lemma}

\begin{lemma}[Lemma 1 of \cite{mania2019certainty}]
\label{lem:strongly_convex_prop}
Let $f_1,f_2$ be two $\mu$-strongly convex twice differentiable functions. Let
\[
x_1=\arg\min_x f_1(x),
\qquad
x_2=\arg\min_x f_2(x).
\]
If $\norm{\nabla f_1(x_2)}\leq\epsilon$, then
\(
\norm{x_1-x_2}\leq\frac{\epsilon}{\mu}.
\)
\end{lemma}

\begin{lemma}[Adapted from Lemma 2 of \cite{mania2019certainty}]
\label{lem:opt_solution_lqr_close}
For $i=1,2$, define
\[
f_i(u;x)
=
\frac{1}{2}u^\top R_i u
+
\frac{1}{2}(A_ix+B_iu)^\top P_i(A_ix+B_iu),
\]
where $R_i$ and $P_i$ are positive definite and $\smin(R_i)\geq1$. Let $K_i$ be the unique matrix such that
\(
\arg\min_u f_i(u;x)=K_ix
\)
for every vector $x$. Let
\[
\Gamma
=
1+\max\{
\norm{A_1},\norm{B_1},\norm{P_1},
\norm{K_1},\norm{R_1^\inv}
\}.
\]
Suppose $0\leq\epsilon<1$ and
\[
\norm{A_1-A_2},
\norm{B_1-B_2},
\norm{R_1-R_2},
\norm{P_1-P_2}
\leq
\epsilon.
\]
Then
\[
\norm{K_1-K_2}
\leq
7\Gamma^4\epsilon.
\]
\end{lemma}

\begin{proof}
The gradient with respect to $u$ is
\[
\nabla f_i(u;x)
=
(B_i^\top P_iB_i+R_i)u+B_i^\top P_iA_ix.
\]
A direct perturbation bound gives
\begin{align*}
\norm{
B_1^\top P_1B_1+R_1
-
B_2^\top P_2B_2-R_2
}
& \leq
7\Gamma^2\epsilon \\
\norm{
B_1^\top P_1A_1
-
B_2^\top P_2A_2
}
& \leq
7\Gamma^2\epsilon.
\end{align*}
Thus, for any $\norm{x}\leq1$,
\[
\norm{\nabla f_1(u;x)-\nabla f_2(u;x)}
\leq
7\Gamma^2\epsilon(\norm{u}+1).
\]
Since $\norm{u_1}\leq\norm{K_1}\norm{x}\leq\norm{K_1}$ and $f_2$ is strongly convex with parameter at least $1$, \cref{lem:strongly_convex_prop} gives
\[
\norm{(K_1-K_2)x}
=
\norm{u_1-u_2}
\leq
7\Gamma^4\epsilon.
\]
Taking the supremum over $\norm{x}\leq1$ proves the result.
\end{proof}

The following lemma is adapted from \cite{mania2019certainty} to additionally consider perturbations in the cost matrices $Q$ and $R$.

\begin{lemma}[Adapted from Proposition 1 of \cite{mania2019certainty}]
\label{lem:control_close_and_tau}
Let $\epsilon>0$ satisfy $\norm{\Delta_A}, \norm{\Delta_B}, \norm{\Delta_R} \leq \epsilon$.
Let
\[
P=\dare(A,B,Q,R),
\qquad
P_\epsilon=\dare(A_\epsilon,B_\epsilon,Q_\epsilon,R_\epsilon),
\]
and suppose that $\norm{P-P_\epsilon}\leq f(\epsilon)$ for some $f(\epsilon)\geq\epsilon$. If $K$ optimizes $\loss(\cdot,Q,R)$ under $(A,B)$ and $K_\epsilon$ optimizes the corresponding LQR objective under $(A_\epsilon,B_\epsilon,Q_\epsilon,R_\epsilon)$, then
\[
\norm{K-K_\epsilon}
\leq
7\Gamma^4 f(\epsilon).
\]
Moreover, if $\rho(A+BK)<\gamma<1$ and
\[
7\Gamma^4 f(\epsilon)
\leq
\frac{1-\gamma}{
2\tau(A+BK,\gamma)\norm{B}
},
\]
then
\[
\tau \bigl(A+BK_\epsilon,\frac{1+\gamma}{2}\bigr)
\leq
\tau(A+BK,\gamma).
\]
\end{lemma}

\begin{proof}
The first claim follows from \cref{lem:opt_solution_lqr_close}, using $\smin(R)\geq1$. For the stability bound, apply \cref{lem:tau_bound_sum} with
\[
M=A+BK,
\qquad
\Delta=B(K_\epsilon-K).
\]
For every $k\geq1$,
\ifdefined\arxiv
\begin{align*}
\norm{(A+BK_\epsilon)^k}
&\leq
\tau(A+BK,\gamma) 
\left(
\tau(A+BK,\gamma)
\norm{B(K_\epsilon-K)}
+\gamma
\right)^k \\
&\leq
\tau(A+BK,\gamma)
\left(\frac{1+\gamma}{2}\right)^k,
\end{align*}
\else
\begin{align*}
\norm{(A+BK_\epsilon)^k}
&\leq
\tau(A+BK,\gamma) \\
&\quad\cdot
\left(
\tau(A+BK,\gamma)
\norm{B(K_\epsilon-K)}
+\gamma
\right)^k \\
&\leq
\tau(A+BK,\gamma)
\left(\frac{1+\gamma}{2}\right)^k,
\end{align*}
\fi
where the last step uses the assumed upper bound on $7\Gamma^4f(\epsilon)$. Multiplying both sides by $((1+\gamma)/2)^{-k}$ and taking the supremum over $k$ proves the claim.
\end{proof}